\documentclass[12pt]{amsart}

\usepackage{amsmath, amssymb, graphics}

\newcommand{\mathsym}[1]{{}}
\newcommand{\unicode}[1]{{}}
\newtheorem{thm}{Theorem}[section]
\newtheorem{cor}[thm]{Corollary}
\newtheorem{lem}[thm]{Lemma}
\newtheorem{prop}[thm]{Proposition}
\theoremstyle{definition}

\newtheorem{defn}[thm]{Definition}
\newtheorem{rem}[thm]{Remark}

\newtheorem*{defn*}{Definition}
\newtheorem*{rems*}{Remarks}
\newtheorem*{rem*}{Remark}

\numberwithin{equation}{section}
\begin{document}

\title [Singularities of GCS of  Lagrangian  submanifolds] {Singularities of equidistants and \\global centre symmetry sets
\\of   Lagrangian  submanifolds}

\author[Domitrz \& Rios]{Wojciech Domitrz \& Pedro de M. Rios}

\address{Warsaw University of Technology, Faculty of Mathematics and Information Science, Plac Politechniki 1, 00-661 Warszawa, Poland}
\email{domitrz@mini.pw.edu.pl}
\address{Departamento de Matem\'atica, ICMC, Universidade de S\~ao Paulo; S\~ao Carlos, SP, 13560-970, Brazil}
\email{prios@icmc.usp.br}

\thanks{W. Domitrz was supported by FAPESP, during his stay in S\~ao Carlos, and by Polish MNiSW grant  N N201 397237. P. de M. Rios received partial support by
FAPESP for a visit to Warsaw.}

\subjclass{57R45, 58K40, 53D12, 58K25, 58K50.}

\keywords{Centre symmetry set, Symplectic geometry, Lagrangian singularities}

\maketitle

\begin{abstract}
We define the global centre symmetry set ($GCS$) of a smooth
closed $m$-dimensional submanifold $M\subset\mathbb R^n$, $n \leq
2m$, which is an affinely invariant generalization of the centre
of a $k$-sphere in $\mathbb R^{k+1}$. The $GCS$ includes both the
centre symmetry set defined by Janeczko \cite{Jan} and the Wigner
caustic defined by Berry \cite{Ber}. We study singularities of the $GCS$ of a Lagrangian submanifold of $\mathbb R^{2m}$ with canonical
symplectic form. The definition of the $GCS$, which slightly
generalizes one by Giblin and Zakalyukin \cite{GZ1}-\cite{GZ3}, is
based on the notion of affine equidistants, so, we first study
singularities of affine equidistants of Lagrangian submanifolds.
Then, we classify
affine-Lagrangian stable singularities of the $GCS$ of  Lagrangian
submanifolds and show that, already for smooth closed convex
curves in $\mathbb R^2$, many singularities of the $GCS$ which are
affine stable are not affine-Lagrangian stable.
\end{abstract}

\section{Introduction}

An $(n-1)$-sphere in $\mathbb R^{n}$ is usually defined as the set of all points in $\mathbb R^{n}$
which are equidistant to a fixed point. Naturally, this point is
called the centre of the sphere, or the centre of
symmetry of the sphere. However, this definition
depends on a choice of metric in $\mathbb R^{n}$.

Now, looking at a circle
on the plane, or even an ellipse, its centre can instead be defined
as the set (in this case consisting of a single element) of
midpoints of straight lines connecting pairs of points on the
curve with parallel tangent vectors.
For a generic smooth convex closed curve, this set is not a single
point, but forms a curve with an odd number of cusps, in the
interior of the smooth original curve. This singular inner curve
has been known as the {\it Wigner caustic} of the  smooth curve since
the work of Berry in the $70$'s,  because of its prominent
appearance in the semiclassical limit of the Wigner function of a
pure quantum state whose classical limit corresponds to the given
smooth curve in $\mathbb R^2$, with canonical symplectic
structure \cite{Ber} \cite{OH}.
Thus, the Wigner caustic  is an affine-invariant
generalization of the centre  of a circle, which extends to higher dimensional smooth closed
submanifolds of $\mathbb R^n$.

On the other hand, the centre of a circle or an ellipse in
$\mathbb R^2$ can also be described as the envelope of all
straight lines connecting pairs of points on the curve with
parallel tangent vectors.
For a generic smooth convex closed curve, this set is not a
single point, but forms a curve with an odd number of cusps, in the
interior of the smooth original curve. This singular  inner curve
has been known as the {\it centre symmetry set} (CSS) of a smooth closed
curve  since the work of Janeczko in the $90$'s
and is an affine-invariant generalization of the centre of
 a circle, which extends to higher
dimensional smooth closed submanifolds of $\mathbb R^n$
\cite{Jan}.

However, except for circles or ellipses, when both symmetry sets are the same point,
the Wigner caustic and the centre symmetry set of a smooth convex closed curve are
not the same singular curve. Instead, the Wigner caustic is
interior to the centre symmetry set and the cusp points of the
inner curve touches the outer one in its smooth part.
A larger centre symmetry set, containing the two previous ones, can be defined an
affine-invariant way, for an arbitrary smooth closed
$m$-dimensional submanifold $M$ of $\mathbb R^n$, for $n\leq 2m$.
We call this new set the {\it global centre symmetry} set of $M$,
denoted by $GCS(M)$.

Our definition is a slight modification of a
definition introduced by Giblin and Zakalyukin
\cite{GZ1}-\cite{GZ3} to study singularities of centre symmetry
sets of hypersurfaces. A key notion in their definition is that of
an affine $\lambda$-equidistant to the smooth submanifold $M$, denoted $E_{\lambda}(M)$ of
which the Wigner caustic is the case $\lambda=1/2$. The
singularities of $E_{\lambda}(M)$ are then fundamental
to characterize $GCS(M)$ and its own
singularities.

In this paper, we
study singularities of
$E_{\lambda}(L)$ and $GCS(L)$, when $L$ a smooth closed {\it Lagrangian}
submanifold of
$(\mathbb{R}^{2n}, \omega)$, where $\omega$ is the canonical symplectic form.
This paper is organized as follows.

In section \ref{defGCS} we
present the definition of affine $\lambda$-equidistants and the global centre symmetry set.
In section \ref{gen-fam-section} we obtain
 the generating families
for the affine equidistants $E_{\lambda}(L)$, relating their general
classification to the well known classification by
 Lagrangian equivalence \cite{AGV}. This is used in section
 \ref{sing-equi-section}
to study  singularities of affine equidistants.
Theorem \ref{equi-stable} states that any caustic of simple stable
 Lagrangian singularity is realizable as $E_{\lambda}(L)$, for some
 Lagrangian $L\subset(\mathbb{R}^{2m},\omega)$.

 In section \ref{criminant-section}
we give a geometric characterization for the criminant of $GCS(L)$
similar to results in \cite{GZ1}-\cite{GZ3} for hypersurfaces. In
section \ref{section-sing-GCS} we introduce the equivalence
relation (also as an equivalence of generating families) that is used to classify
the singularities of $GCS(L)$. We show that only
singularities of the criminant, the smooth part of the Wigner
caustic, or tangent union of both, are affine-Lagrangian stable.

Finally, section \ref{css-class-Lag-curves} is devoted to the
study of the GCS of Lagrangian curves. First, we state two
theorems for the GCS of convex curves in $\mathbb R^2$ when no
symplectic structure is considered. The results presented in
Theorem \ref{curvestability} are not new (\cite{Ber}, \cite{Jan},
\cite{GH}-\cite{Gib2}), but, in Theorem \ref{globalresults} the inequality on the
number of cusps of the CSS and the Wigner caustic had not been
mentioned before. Pictures illustrate these theorems. Then, we
 show that most of the singularities
which were affine-stable when no symplectic structure was
considered are not affine-Lagrangian stable.
In other words, although any smooth curve on $\mathbb R^2$ is Lagrangian,
the singularities of their GCS are sensitive to the presence
of a symplectic form to be accounted for, that is, there is a breakdown
of the stability of these singularities. This is  similar to some results  in
\cite{D}-\cite{DR}.

\

\noindent{\bf Acknowledgements:} We specially thank M.A.S. Ruas
for many stimulating discussions and invaluable remarks. We also thank P. Giblin and S.
Janeczko for discussions and V. Goryunov
for remarks.

\section{Definition of the global centre symmetry
set.}\label{defGCS}

Let $M$ be a smooth closed $m$-dimensional submanifold of
the affine space $\mathbb R^{n}$, with $n\leq 2m$. Let $a, b$ be points of $M$.
Let $\tau_{a-b}$ be the
translation by the vector $(a-b)$, i.e.,
$\tau_{a-b}:\mathbb R^n \ni x\mapsto x+(a-b) \in \mathbb R^n.$
\begin{defn}\label{parallelism} A pair  $a, b \in M$ ($a\ne b$) is  a
{\bf weakly parallel} pair if
$$T_aM + \tau_{a-b}(T_bM)\ne T_a\mathbb R^n.$$
A weakly parallel pair $a, b \in M$ is
called {\bf $k$-parallel} if
$$\dim(T_aM \cap \tau_{b-a}(T_bM))=k.$$
If $k=m$ the pair $a, b \in M$ is called {\bf strongly parallel},
or just {\bf parallel}. We also refer to $k$ as the {\bf degree of
parallelism} of the pair $(a,b)$. \end{defn}

\begin{defn}
A {\bf chord} passing through a pair $a,b$, is the line
$$
l(a,b)=\{x\in \mathbb R^n|x=\lambda a + (1-\lambda) b, \lambda \in
\mathbb R\}.
$$
\end{defn}
\begin{defn} For a given $\lambda$, an {\bf affine
$\lambda$-equidistant} of $M$, $E_{\lambda}(M)$, is the set of all
$x\in \mathbb R^n$ such that $x=\lambda a + (1-\lambda) b$, for
all weakly parallel pairs $a,b \in M$. $E_{\lambda}(M)$ is also
called a (affine) {\bf momentary equidistant} of $M$.
Whenever $M$ is understood, we write
$E_{\lambda}$ for $E_{\lambda}(M)$.
\end{defn}
Note that, for any
$\lambda$, $E_{\lambda}(M)=E_{1-\lambda}(M)$ and in particular
$E_0(M)=E_1(M)=M$. Thus, the case $\lambda=1/2$ is special:
\begin{defn}$E_{{1}/{2}}(M)$ is called the {\bf Wigner caustic} of $M$ \cite{Ber} \cite{OH}.
\end{defn}
The {\it extended affine space} is the space $\mathbb
R^{n+1}_e=\mathbb R\times \mathbb R^{n}$ with coordinate
$\lambda\in \mathbb R$ (called {\it affine
time}) on the first factor and projection on the second factor denoted by
$\pi:\mathbb R^{n+1}_e\ni (\lambda,x)\mapsto
x \in \mathbb R^{n}$.
\begin{defn}\label{aewf}  The {\bf affine extended wave front} of $M$, $\mathbb E(M)$,
is the union of
all affine equidistants each embedded into its own slice of the
extended affine space: $\mathbb E(M)=\bigcup_{\lambda\in \mathbb R}
\ \{\lambda\}\times E_{\lambda}(M) \ \subset \mathbb R_e^{n+1}.$
\end{defn}

Note that, when $M$ is a circle on the plane, $\mathbb E(M)$
is the (double) cone, which is a smooth manifold with nonsingular projection
$\pi$ everywhere, but at its singular point, which projects to the centre
of the circle. From this, we generalize the notion of centre of symmetry.
Thus, let $\pi_r$ be the restriction of $\pi$ to the affine extended
wave front of $M$: $\pi_r=\pi|_{\mathbb E(M)}$. A point $x\in \mathbb E(M)$ is a
{\bf critical} point of $\pi_r$ if the germ of $\pi_r$ at $x$
fails to be the germ of a regular projection of a smooth
submanifold. We now introduce the main definition of this paper:
\begin{defn}\label{GCS}\label{GCSparts}
The {\bf global centre symmetry} set of $M$,  $GCS(M)$, is the
image under $\pi$  of the locus of critical points of $\pi_r$.
\end{defn}
\begin{rem} The set $GCS(M)$ is the bifurcation set of a family
of affine equidistants (family of chords of weakly parallel pairs)
of $M$.

In general, $GCS(M)$ consists of two components: the {\bf caustic
$\Sigma(M)$} being the projection of the singular locus of $\mathbb E(M)$
and the {\bf criminant $\Delta(M)$} being the (closure of) the
image under $\pi_r$ of the set of regular points of $\mathbb E(M)$ which
are critical points of the projection $\pi$ restricted to the
regular part of $\mathbb E(M)$. $\Delta(M)$ is the envelope of the family
of regular parts of momentary equidistants, while $\Sigma(M)$
contains all the singular points of momentary equidistants.
\end{rem}

The above definition (with its following remarks) is a
slight modification of the definition that has already been
introduced by Giblin and Zakalyukin \cite{GZ1}. However, in our present
definition the whole manifold $M$ is considered, as opposed to
pairs of germs, as in \cite{GZ1}, and weak parallelism is also
taken into account. Considering the
whole manifold in the definition leads to the following simple but
important result:

\begin{thm}\label{WinGCS}
The set $GCS(M)$ contains the Wigner caustic of $M$.
\end{thm}
\begin{proof}
Let $x$ be a regular point of $E_{\frac{1}{2}}(M)$. Then
$x=\frac{1}{2}(a+b)$ for a weakly parallel pair $a,b \in M$. It
means that $x$ is a intersection point of the chords $l(a,b)$ and
$l(b,a)$. Then $\mathbb E(M)$ contains the sets
$$
\{(\lambda,\lambda a+(1-\lambda)b)|\lambda \in \mathbb R\}, \
\{(\lambda,(1-\lambda) a+\lambda b)|\lambda \in \mathbb R\}.
$$
If $(\frac{1}{2},x)$ is a regular point of $\mathbb E(M)$ then the above
sets are included in the tangent space to $\mathbb E(M)$ at
$(\frac{1}{2},x)$. Therefore a fiber $\{(\lambda,x)|\lambda
\in \mathbb R\}$ is included in the tangent space of $\mathbb E(M)$. Thus
if $(\frac{1}{2},x)$ is a regular point of $\mathbb E(M)$ then $x$ is in
the criminant $\Delta(M)$. If $(\frac{1}{2},x)$ is not a regular
point of $\mathbb E(M)$  then $x$ is in the caustic $\Sigma(M)$.
\end{proof}

 If  $M\subset\mathbb R^{2}$ is a smooth curve, $E_{1/2}(M)\ni x$ is the bifurcation set for the number of chords connecting two points in $M$ given a
chord midpoint $x\in\mathbb R^{2}$  \cite{Ber}. Similarly, if ${\mathcal R}_x:\mathbb R^{2}\to\mathbb R^{2}$
denotes reflection through $x\in\mathbb R^{2}$, then $x\in E_{1/2}(M)$ when $M$ and  ${\mathcal R}_x(M)$ are not
transversal \cite{OH}\cite{GJ}. Finally, let $A(x,\kappa)$  be  the area of the planar region enclosed by $M$ and a chord as a function of a point $x$ on the chord and a variable $\kappa$ locating one of the endpoints of the chord on the curve. Then,
 $A(x,\kappa)$ is a generating family for $E_{1/2}(M)$  \cite{Ber,Gib2}, which
 is generalized below  to every $\lambda$-equidistant of any Lagrangian submanifold.

\section{Generating families}\label{gen-fam-section}

Consider the product affine space: $\mathbb R^{n}\times \mathbb R^{n}$ with coordinates $(x_+,x_-)$
and the tangent bundle to $\mathbb R^{n}$:  $T\mathbb
R^{n}=\mathbb R^{n}\times \mathbb R^{n}$ with coordinate system
$(x,\dot{x})$ and standard projection
$pr: T\mathbb R^{n}\ni (x,\dot{x})\rightarrow x \in \mathbb R^{n}$.

\begin{defn} $\forall \lambda\in\mathbb{R}\setminus  \{0,1\}$, a {\bf $\lambda$-chord transformation}
$$\Phi_{\lambda}:\mathbb{R}^{n}\times\mathbb{R}^{n}\to T\mathbb{R}^{n} \ , \ (x^+,x^-)\mapsto(x,\dot{x})$$
is a {\it linear} diffeomorphism  defined by
the {\it $\lambda$-point equation}:
\begin{equation}\label{x}
x=\lambda x^+ + (1-\lambda)x^- \ ,
\end{equation}
for the $\lambda$-point $x$, and a {\it chord equation}:
\begin{equation}\label{genx.}
\dot{x}=\lambda x^+ - (1-\lambda) x^-.
\end{equation}
\end{defn}

Now, let $M$ be a smooth closed $m$-dimensional submanifold of the
affine space $\mathbb R^{n}$ ($2m\ge n$) and consider the product
$M\times M\subset \mathbb R^{n}\times\mathbb R^{n}$. Let $\mathcal
M_{\lambda}$ denote the image of $M\times M$ by
a  $\lambda$-chord transformation,
$$\mathcal M_{\lambda} = \Phi_{\lambda}(M\times M) \ .$$
\begin{thm}\label{gensing}
The set of critical values of the standard projection $pr:
T\mathbb R^{n}\to\mathbb R^{n}$ restricted to $\mathcal
M_{\lambda}$ is $E_{\lambda}(M)$.
\end{thm}
\begin{proof} If $a$ is a critical value of $pr|_{\mathcal M_{\lambda}}$, $\dim T_{(a,\dot
a)}\mathcal M_{\lambda}\cap T_{(a,\dot
a)}pr^{-1}(a)$ is greater than $2m-n$. Let $v_1,\cdots, v_k$ for $k>2m-n$ be a basis of $
T_{(a,\dot{a})}\mathcal M_{\lambda}\cap
T_{(a,\dot{a})}pr^{-1}(a)$, of the
form $v_j=\sum_{i=1}^{n}\alpha_{ji} \frac{\partial}{\partial
\dot{x}_i}|_{(a,\dot{a})}$ for $j=1,\cdots,k$ . But
we have that
$(\Phi_{\lambda}^{-1})_{\ast}(v_j)=\frac{1}{2\lambda}v^+_j-\frac{1}{2(1-\lambda)}v^-_j$,
where
$$v^+_j= \sum_{i=1}^{n}\alpha_{ji}
\frac{\partial}{\partial x^+_i}|_{a^+}\in T_{a^+}M , \ \
v^-_j=
\sum_{i=1}^{n}\alpha_{ji} \frac{\partial}{\partial
x^-_i}|_{a^-}\in T_{a^-}M.$$
It implies that  $v_j^+\in
T_{a^+}M\cap \tau_{(a^+-a^-)}T_{a^-}M$ for $j=1,\cdots,k$. Since $k>2m-n$ then
$T_{a^+}M+\tau_{(a^+-a^-)}T_{a^-}M\ne T_{a^+}\mathbb R^n$ and
consequently $a^+, a^-$ is a $k$-parallel pair.
Hence $\lambda a^++(1-\lambda)a^-=a\in E_{\lambda}$.

Now, assume $a\in E_{\lambda}$. Then $a=\lambda
a^++(1-\lambda)a^-$ for a weakly $k$-parallel pair $a^+, a^-$ for $k>2m-n$.
Thus there exist linearly independent vectors
$v^+_j=\sum_{i=1}^{n}\alpha_{ji} \frac{\partial}{\partial
x^+_i}|_{a^+}\in T_{a^+}M\cap \tau_{(a^+-a^-)}T_{a^-}M$ for
$j=1,\cdots,k$. Consider linearly independent  vectors
$v_j=(\Phi_{\lambda})_{\ast}((1-\lambda)v^+_j-
\lambda\tau_{(a^--a^+)}v^+_j)$ for $j=1,\cdots,k$. Then, $v_j$ belongs to $T_{(a,\dot{a})}\mathcal
M_{\lambda}$ and $pr_{\ast}(v_j)=0$ for
$j=1,\dots, k$. Thus $a$ is a critical value of $pr|_{\mathcal
M_{\lambda}}$.
\end{proof}

Let $(\mathbb R^{2m}, \omega)$ be the affine symplectic space with
canonical coordinates ${p_i,q_i}$, so that
$\omega=\sum_{i=1}^m dp_i\wedge dq_i$, and let $L$ be  a smooth
closed  Lagrangian submanifold of $(\mathbb R^{2m}, \omega)$. For a fixed $\lambda\in \mathbb R\setminus  \{0,1\}$, consider
the product affine space $\mathbb R^{2m}\times \mathbb R^{2m}$
with the $\lambda$-weighted symplectic form
\begin{equation}\label{symplform}
\delta_{\lambda}\omega=2\lambda^2
\pi_1^{\ast}\omega-2(1-\lambda)^2\pi_2^{\ast}\omega \ ,
\end{equation}
where $\pi_i$ is the projection of $\mathbb R^{2m}\times \mathbb R^{2m}$
on $i$-th factor for $i=1,2$.

Now, let $\Phi_{\lambda}$ be the $\lambda$-chord
transformation (\ref{x})(\ref{genx.}). Then,
\begin{equation}\label{Omega'}
\left(\Phi_{\lambda}^{-1}\right)^{\ast}(\delta_{\lambda}\omega)
 \ = \ \dot{\omega} \ .
\end{equation}
where $\dot{\omega}$
is the canonical symplectic form on the tangent bundle to $(\mathbb{R}^{2m},\omega)$,
defined by $\dot{\omega}(x,\dot{x})=d\{\dot{x}\lrcorner\omega\}(x)$ or, in Darboux coordinates,
\begin{equation}\label{omega.}
\dot{\omega}=\sum_{i=1}^m
d\dot{p_i}\wedge dq_i+ dp_i\wedge d\dot{q_i} \ .
\end{equation}

The fibers of $T\mathbb{R}^{2m}$ are  Lagrangian for
$\dot{\omega}$, so that $pr:T\mathbb R^{2m}\rightarrow \mathbb R^{2m}$ defines
a {\bf  Lagrangian fiber bundle} with respect to
$\dot{\omega}$, that is, a fiber bundle whose fibers are Lagrangian in the
total symplectic space.

\begin{prop} The restriction of the projection $pr:T\mathbb R^{2m}\rightarrow \mathbb R^{2m}$ of
$(T\mathbb{R}^{2m},\dot{\omega})$ to the Lagrangian submanifold
$$\mathcal L_{\lambda}=\Phi_{\lambda}(L\times L)$$
is a  Lagrangian map \cite{AGV}.  The set of critical values of a  Lagrangian map is called a {\bf
caustic} and (Theorem \ref{gensing})  the caustic of $pr|_{\mathcal L_{\lambda}}$ is
$E_{\lambda}(L)$.
\end{prop}

\begin{defn} $E_{\lambda}(L)$ and $E_{\lambda}(\widetilde
L)$ are {\bf  Lagrangian equivalent} if the  Lagrangian maps
$pr|_{\mathcal L_{\lambda}}$ and $pr|_{\widetilde{\mathcal
L}_{\lambda}}$ are  Lagrangian equivalent (see \cite{AGV}).
\end{defn}

It follows from above definitions:

\begin{prop}\label{affine-symplectic-inv-Lag-equiv}
The classification of $E_{\lambda}(L)$ by  Lagrangian equivalence
is affine symplectic invariant, i.e., invariant under the standard action of the affine symplectic group on $(\mathbb R^{2m},\omega)$.
\end{prop}

From the above, we also use the term  {\bf
affine-Lagrangian equivalence}  for Lagrangian equivalence
(see \cite{AGV}) of $E_{\lambda}(L)$.

Now, let $L^+$ and $L^-$ denote germs of $L$ at points $a^+$ and $a^-$.

\begin{prop} \label{s}
If the pair $a^+,a^-$ is
$k$-parallel, there exists canonical
coordinates $(p,q)$ on $\mathbb R^{2m}$ and function germs $S^+$
and $S^-$ such that
\begin{equation}\label{s-}
L^+: p_i=\frac{\partial S^+}{\partial q_i}(q_1,\cdots,q_m), \
 \  i=1,\cdots,m
\end{equation}
\begin{equation}
L^-: \begin{cases}p_j= \ \ \frac{\partial S^-}{\partial
q_j}(q_1,\cdots,q_k,p_{k+1},\cdots,p_m),    \  \ j=1,\cdots,k, \\
q_l=-\frac{\partial S^-}{\partial
p_l}(q_1,\cdots,q_k,p_{k+1},\cdots,p_{m}),  \  \
l=k+1,\cdots,m
\end{cases}\nonumber
\end{equation}
and $d^2S^+(q_{a,1}^{+},\cdots,q_{a,m}^+)=0$ and
$d^2S^-(p_{a,1}^-,\cdots,p_{a,k}^-,q_{a,k+1}^-,\cdots,p_{a,m}^-)=0$,
where $a^+=(p_a^+,q_a^+)$ and $a^-=(p_a^-,q_a^-)$.
\end{prop}

\begin{proof}
We can find a linear symplectic change of coordinates such that
$T_{a^+}L^+=\{p=p^{+}_a\}$, where $a^+=(p_a^+,q_a^+)$, and
$T_{a^-}L^-=\{p_1=p^{-}_{a,1},\cdots,p_k=p^-_{a,k},q_{k+1}=q^-_{a,k+1},\cdots,q_m=q^-_{a,m}\}$,
where  $a^-=(p_a^-,q_a^-)$. Since $L$ is a smooth
 Lagrangian submanifold, it follows from standard considerations
that it can be described locally by differentials of generating
functions of the forms stated above in neighborhoods of $a^+$ and
$a^-$, in which case we have that $d^2S^+|a^+=d^2S^-|a^-=0$.
\end{proof}

Let the arguments
of the function $S^+$ be denoted by $(q_1^+,\cdots,q_m^+)$
and the arguments of the function $S^-$ by
$(q_1^-,\cdots,q_k^-,p^-_{k+1},\cdots,p^-_m)$.
Let $q=(q_1,\cdots,q_m)$, $p=(p_{1},\cdots,p_m)$, $\dot q=(\dot
q_1,\cdots,\dot q_m)$, $\dot p=(\dot p_{1},\cdots,\dot p_m)$.

Also, let $\beta=(\beta_1,\cdots,\beta_m)$
and, for any $k<m$, let $[k]=\{1,\cdots,k\}$, so that $\beta_{[k]}=(\beta_1,\cdots,\beta_k)$,
and $\alpha_{[m]\setminus [k]}=(\alpha_{k+1},\cdots,\alpha_m)$.

Let $L^+\times L^-$ denote the germ of $L\times L$ at the point
$(a^+,a^-)\in L\times L$ so that $\mathcal
L_{\lambda}=\Phi_{\lambda}(L^+\times L^-)$ is the germ at
$(a,\dot{a})$, where $a=\lambda a^++(1-\lambda)a^-$, $\dot
a=\lambda a^+-(1-\lambda)a^-$, of a smooth  Lagrangian submanifold
of $(T\mathbb
R^{2m}, \dot{\omega})$.

\begin{thm} \label{gen-fam}
If the pair $a^+,a^-$ is
$k$-parallel and germs $L^+$ and $L^-$ are given by
(\ref{s-}) then the germ of the
{\bf generating family}
\begin{eqnarray}\label{genfam}
&F_{\lambda}(p,q,\alpha_{[m]\setminus  {[k]}},\beta)=&\\
 &2\lambda^2
S^+\left(\frac{q+\beta}{2\lambda}\right)-2(1-\lambda)^2
S^-\left(\frac{q_{[k]}-\beta_{[k]} \ , \ p_{[m]\setminus  {[k]}}-\alpha_{[m]\setminus  {[k]}}}{2(1-\lambda)}\right)&  \nonumber\\
&-\sum_{i=1}^k p_i \beta_i +\frac{1}{2}\sum_{j=k+1}^m q_j \alpha_j
-p_j \beta_j - \alpha_j \beta_j -p_j q_j &\nonumber
\end{eqnarray}
generates the germ of $\mathcal L_{\lambda}$ at $(a,\dot{a})$ as follows:
$$
\mathcal L_{\lambda}=\left\{(\dot{p},\dot{q},p,q):\exists
(\alpha,\beta) \ \dot{p}=\frac{\partial F_{\lambda}}{\partial q},
\ \dot{q}=-\frac{\partial F_{\lambda}}{\partial p}, \ \frac{\partial
F_{\lambda}}{\partial \alpha}=\frac{\partial F_{\lambda}}{\partial
\beta}=0 \right\}.
$$
\end{thm}

\begin{proof} The proof is a straightforward calculation.
\end{proof}

\begin{rem} \label{corank}
Note from
(\ref{genfam}) that {\it the degree of parallelism is the corank
of the singularity}, i.e. the corank of the Hessian of the
function
$$\mathbb R^{2m-k}\ni(\alpha_{[m]\setminus  {[k]}},\beta)\mapsto
F_{\lambda}(p_a,q_a,\alpha_{[m]\setminus  {[k]}},\beta)\in \mathbb
R$$
\end{rem}

\begin{thm}[\cite{AGV}]\label{LagstabR+}
Germs of  Lagrangian maps are  Lagrangian equivalent iff the germs of their generating families are stably
$\mathcal R^+$-equivalent.
\end{thm}

\begin{cor}\label{Lagstabfam}
 Germs
$E_{\lambda}(L)$ and $E_{\lambda}(\tilde L)$ are
 Lagrangian equivalent iff germs of
generating families for ${\mathcal L}_{\lambda}$ and $\tilde
{\mathcal L}_{\lambda}$ are stably $\mathcal R^+$-equivalent.
\end{cor}

\section{Singularities of equidistants of  Lagrangian
submanifolds}\label{sing-equi-section}

We have the following results on singularities of affine
equidistants of closed  Lagrangian submanifolds, up to
 Lagrangian equivalence:
\begin{thm}\label{equi-stable} Any caustic of simple stable Lagrangian singularity (A-D-E singularities) in the $4m$-dimensional
symplectic tangent bundle $(T\mathbb R^{2m}, \dot{\omega})$ is
realizable as $E_{\lambda}(L)$, for some smooth closed
 Lagrangian submanifold $L$ in $(\mathbb R^{2m},\omega)$.
\end{thm}

The generic Lagrangian maps for manifolds of dimension smaller than $6$ have only simple stable Lagrangian singularities (\cite{AGV}). Therefore we obtain the following corollary.

\begin{cor}
Any germ of generic caustics on  $2m$-dimensional manifold for $m=1,2$ is realizable as $E_{\lambda}(L)$, for some smooth Lagrangian submanifold $L$ in $(\mathbb R^{2m},\omega)$.
\end{cor}

\begin{proof}[Proof of Theorem \ref{equi-stable}] We use the method described in \cite{AGV}. For a fixed $\lambda$, let $x=(p,q)$ and
$\kappa=(\alpha,\beta)$. From (\ref{genfam}) we easily see that
$$
\text{rank}_{(a,\dot{a})}\left[\frac{\partial^2
F_{\lambda}}{\partial \kappa^2}, \ \frac{\partial^2 F_{\lambda}}{
\partial \kappa \partial x}\right]=2m-k,$$ hence is equal to the
dimension of $\kappa$-space. We find $S^+$ and $S^-$ such that $F_{\lambda}(x,\kappa)$ is a $\mathcal
R^+$-versal deformation of A-D-E singularities. Let
$$S^{+}(q^+)=\sum_{i=1}^m
p_{a,i}^{+}(q_i^{+}-q_{a,i}^{+})+S^{+}_3(q^+-q_a^+)$$
\begin{eqnarray*}
S^{-}(q^-_{[k]},p^-_{[m]\setminus [k]})&=&\sum_{i=1}^k
p_{a,i}^{-}(q_i^{-}-q_{a,i}^{-})-\sum_{i=k+1}^m
q_{a,i}^{-}(p_i^{-}-p_{a,i}^{-})+\\
&+&S^{-}_3(q^-_{[k]}-q^-_{a,[k]}, p^-_{[m]\setminus [k]}-p_{a,[m]\setminus [k]}^-),
\end{eqnarray*}
 where we used Proposition \ref{s} and where $S^{\pm}_3\in \mathfrak m^3$ ($\mathfrak m$ is the maximal ideal of the ring of smooth function-germs on $\mathbb R^n$ at $0$).
We write the generating families in coordinates $\tilde p=p-p_a$,
$\tilde q=q-q_a$, $s=\alpha-\dot p_{a}$, $t=\beta-\dot q_a$, where
$a=(p_a,q_a)$, $\dot{a}=(\dot{p}_a,\dot q_a)$. By Theorem
\ref{gen-fam} we obtain
\begin{eqnarray}\label{F}
&F_{\lambda}(\tilde p,\tilde q,s,t)=&\\
 &2\lambda^2
S^+_3\left(\frac{\tilde q+t}{2\lambda}\right)-2(1-\lambda)^2
S^-_3\left(\frac{\tilde{q}_{[k]}-t_{[k]} \ , \ \tilde{p}_{[m]\setminus [k]}-s_{[m]\setminus [k]}}{2(1-\lambda)}\right)&  \nonumber\\
&-\sum_{i=1}^k \tilde{p}_i t_i +\frac{1}{2}\sum_{j=k+1}^m \tilde{q}_j s_j -\tilde{p}_j t_j - s_j t_j -\tilde{p}_j
\tilde{q}_j  + \sum_{l=1}^m \dot{p}_{a,l}\tilde{q}_l-\dot{q}_{a,l}\tilde{p}_l&\nonumber
\end{eqnarray}
\begin{eqnarray}\label{f}
&f_{\lambda}(s,t)=F_{\lambda}(0,0,s,t)=&\\
&2\lambda^2
S^+_3\left(\frac{t}{2\lambda}\right)-2(1-\lambda)^2
S^-_3\left(\frac{-t_{[k]}, -s_{[m]\setminus [k]}}{2(1-\lambda)}\right) -\frac{1}{2}\sum_{j=k+1}^m s_jt_j&  \nonumber
\end{eqnarray}

\noindent The following singularities are realizable by
generating function-germs:
$$A_{2l}:$$
$$S^{+}_3(\tilde q^+)=\lambda (\tilde q_1^+)^3 + (\tilde
q_1^+)^{2l+1}+\sum_{i=2}^l \tilde q_i^+(\tilde q_1^+)^{2i-1},$$
$$S^{-}_3(\tilde q^-_1,\tilde p^-_2,\cdots,\tilde p^-_m)=-(1-\lambda) (\tilde q_1^-)^3 +\sum_{i=2}^{l-1} \tilde p_i^-(\tilde q_1^-)^{2(l-i+1)}.$$
$$A_{2l+1}:$$
$$S^{+}_3(\tilde q^+)=\lambda (\tilde q_1^+)^3 + (\tilde
q_1^+)^{2l+2}+\sum_{i=2}^l \tilde q_i^+(\tilde q_1^+)^{2i-1},$$
$$S^{-}_3(\tilde q^-_1,\tilde p^-_2,\cdots,\tilde p^-_m)=-(1-\lambda) (\tilde q_1^-)^3 +\sum_{i=2}^{l} \tilde p_i^-(\tilde q_1^-)^{2(l-i+2)}.$$$$D_{2l}:$$
$$S^{+}_3(\tilde q^+)=\lambda (\tilde q_1^+)^3 +\tilde q_2^+ (\tilde q_1^+)^2  \pm (\tilde
q_2^+)^{2l-1}+\lambda (\tilde q_2^+)^3+\sum_{i=2}^{l-1} \tilde
q_{i+1}^+(\tilde q_2^+)^{2i-1},$$
$$S^{-}_3(\tilde q^-_{[2]},\tilde p^-_{[m]\setminus [2]})=-(1-\lambda)(\tilde q_1^-)^3
-(1-\lambda) (\tilde q_2^-)^3+\sum_{i=2}^{l-2} \tilde
p_{i+1}^-(\tilde q_2^-)^{2(l-i)}.$$
$$D_{2l+1}:$$
$$S^{+}_3(\tilde q^+)=\lambda (\tilde q_1^+)^3 +\tilde q_2^+ (\tilde q_1^+)^2  \pm (\tilde
q_2^+)^{2l}+\lambda (\tilde q_2^+)^3+\sum_{i=2}^{l-1} \tilde
q_{i+1}^+(\tilde q_2^+)^{2i-1},$$
$$S^{-}_3(\tilde q^-_{[2]},\tilde p^-_{[m]\setminus [2]})=-(1-\lambda)(\tilde q_1^-)^3
-(1-\lambda) (\tilde q_2^-)^3+\sum_{i=2}^{l-1} \tilde
p_{i+1}^-(\tilde q_2^-)^{2(l-i+1)}.$$
$$E_6:$$
$$S^{+}_3(\tilde q^+)=(\tilde q_1^+)^3 \pm (\tilde
q_2^+)^{4}+ \lambda\tilde q_1^+ (\tilde q_2^+)^2  +\lambda (\tilde
q_2^+)^3+ \tilde q_{1}^+(\tilde q_2^+)^{2}\tilde q_3^+,$$
$$S^{-}_3(\tilde q^-_{[2]},\tilde p^-_{[m]\setminus [2]})=-(1-\lambda)\tilde q_1^-(\tilde
q_2^-)^2 -(1-\lambda) (\tilde q_2^-)^3.$$
$$E_7:$$
$$S^{+}_3(\tilde q^+)=(\tilde q_1^+)^3 +\tilde q_1^+(\tilde
q_2^+)^{2}+ \lambda\tilde q_1^+ (\tilde q_2^+)^2  +\lambda (\tilde
q_2^+)^3+ (\tilde q_2^+)^{3}\tilde q_3^+,$$
$$S^{-}_3(\tilde q^-_{[2]},\tilde p^-_{[m]\setminus [2]})=-(1-\lambda)\tilde q_1^-(\tilde
q_2^-)^2 -(1-\lambda) (\tilde q_2^-)^3+(\tilde q_2^-)^{4}\tilde
p_3^-.$$
$$E_8:$$
$$S^{+}_3(\tilde q^+)=(\tilde q_1^+)^3 +(\tilde
q_2^+)^{5}+ \lambda\tilde q_1^+ (\tilde q_2^+)^2  +\lambda (\tilde
q_2^+)^3+ \tilde q_1^+(\tilde q_2^+)^{2}\tilde q_3^++\tilde
q_1^+(\tilde q_2^+)^{3}\tilde q_4^+,$$
$$S^{-}_3(\tilde q^-_{[2]},\tilde p^-_{[m]\setminus [2]})=-(1-\lambda)\tilde q_1^-(\tilde
q_2^-)^2 -(1-\lambda) (\tilde q_2^-)^3+(\tilde q_2^-)^{3}\tilde
p_3^-.$$

By long but straightforward calculations one can show that (\ref{F}) is a $\mathcal R^+$-versal deformation of (\ref{f}) for the above choices of $S^{\pm}_3$.
\end{proof}

\section{The GCS of a Lagrangian submanifold: the
criminant}\label{criminant-section}

We now begin the study of singularities of the global centre symmetry set of a smooth closed  Lagrangian submanifold
$L\subset(\mathbb{R}^{2m},\omega)$.

The classification of the Wigner caustic of $L$ has been mostly carried out in the last
section. In a subsequent paper \cite{DRs2}, we study
$E_{1/2}(L)$ in a neighborhood of $L$, considering pairs of points of the type $(a,a)\in L\times L$
as strongly parallel pairs. In terms of the generating families of section 4, these are odd functions of the variables, so
we must consider classification {\it in the category of odd functions} \cite{DRs2}.
It also implies a hidden $\mathbb Z_2$-symmetry for singularities for the Wigner caustic on shell.

In order to study the global centre symmetry set,  the whole $\lambda$-family must be considered together.
Due to the Lagrangian condition, we resort to a classification via
generating families.
We
know that $E_{\lambda}(L)$ is the caustic of  $\mathcal L_{\lambda}=\Phi_{\lambda}(L\times L)$. The generating family for $\mathcal
L_{\lambda}$ is given by $F_{\lambda}(p,q,\alpha,\beta)$ of the
form (\ref{genfam}). Since $\mathbb E(L)$ is the union of $\{\lambda\}\times
E_{\lambda}$, the germ of $\mathbb E(L)$ is
described in the following way (for $\kappa=(\alpha,\beta)$):
\begin{prop}\label{wl}
$\mathbb E(L)=\left\{(\lambda,p,q):\exists \kappa \ \frac{\partial
F_{\lambda}}{\partial \kappa}=0, \  \det\left[\frac{\partial^2
F_{\lambda}}{\partial \kappa_i \partial
\kappa_j}\right]=0\right\}$.
\end{prop}
Let us consider the fiber bundle
\begin{equation}\label{Projection}
Pr:T^{\ast}\mathbb R \times
T\mathbb R^{2m}\ni
((\lambda^*,\lambda),(\dot{p},\dot{q},p,q))\mapsto
(\lambda,(p,q))\in \mathbb R\times \mathbb R^m.
\end{equation}

The above bundle with the canonical symplectic structure
$$d\lambda^*\wedge d\lambda+\dot{\omega}$$ is a Lagrangian
fiber bundle. For $F_{\lambda}$ given by (\ref{genfam}) in Theorem
\ref{gen-fam}, let
\begin{equation}\label{FF} F(\lambda,p,q,\alpha,\beta)=F_{\lambda}(p,q,\alpha,\beta).\end{equation}
\begin{prop}\label{mathcalL}
The germ of $\mathbb E(L)$ is the caustic of the germ of a
Lagrangian submanifold $\mathcal L$ of
$(T^{\ast}\mathbb R \times T\mathbb R^{2m},d\lambda^*\wedge d\lambda+\dot\omega)$ generated by the family
$F$ given by (\ref{genfam})-(\ref{FF}) in the following way ($\kappa=(\alpha,\beta)$):
\begin{equation}\label{mathcal-L}
\left\{((\lambda^*,\lambda),(\dot{p},\dot{q},p,q)): \exists
\kappa \  \ \lambda^*=\frac{\partial F}{\partial \lambda}, \
\dot{p}=\frac{\partial F}{\partial q}, \ \dot{q}=-\frac{\partial
F}{\partial p}, \ \frac{\partial F}{\partial \kappa}=0 \right\}.
\end{equation}
\end{prop}

\subsection{Geometric characterization of the criminant of the GCS of a Lagrangian submanifold}
Remind that the criminant $\Delta(L)$ is the (closure of) the image under
$\pi_r$ of the set of regular points of $\mathbb E(L)$ which are
critical points of the projection $\pi$ restricted to the regular
part of $\mathbb E(L)$. That is, the criminant $\Delta(L)$ is the envelope
of the family of regular parts of momentary equidistants. We find
the condition for the tangency to the fibers of the projection
$\pi:(\lambda,p,q)\mapsto (p,q)$.
\begin{prop}\label{regular}
If  $(\lambda,a)$ is a regular point of $\mathbb E(L)$ then there
exists a $1$-parallel pair $a^+, a^-$ such that $a=\lambda a^+
+(1-\lambda) a^-$.
\end{prop}
\begin{proof}
If $(\lambda,a)$ is a regular point of $\mathbb E(L)$ then
the rank of the map
\begin{equation}\label{maprank}
\kappa \mapsto
\left(\frac{\partial F}{\partial
\kappa}(\lambda_a,p_a,q_a,\kappa), \ \det\left[\frac{\partial^2
F}{\partial \kappa_i
\partial \kappa_j}(\lambda_a,p_a,q_a,\kappa)\right]\right)
\end{equation}
is maximal $2m-k$. It implies that
$\text{corank}\left[\frac{\partial^2 F}{\partial \kappa_i
\partial \kappa_j}(\lambda_a,p_a,q_a,\kappa_a)\right]$ is $1$. By
Remark \ref{corank} we obtain that $a^+, a^-$ is a $1$-parallel
pair.
\end{proof}
\begin{prop}\label{tangent}
Let $(\lambda_a,a)=(\lambda_a,p_a,q_a)$ be a regular point of
$\mathbb E(L)$. The fiber of $\pi_r=\pi|_{\mathbb E(M)}$ is tangent to $\mathbb
E(L)$ at $(\lambda_a,a)$ if and only if
\begin{equation}\label{rank} \text{rank}\left[\frac{\partial^2
F}{\partial \lambda
\partial \kappa_j}, \frac{\partial^2 F}{\partial \kappa_i
\partial \kappa_j}\right]=\text{rank}\left[\frac{\partial^2
F}{\partial \kappa_i \partial \kappa_j}\right]=2m-2
\end{equation}
at $(\lambda_a,p_a,q_a,\kappa_a)$ s.t. $\frac{\partial
F}{\partial \kappa}(\lambda_a,p_a,q_a,\kappa_a)=\det\left[\frac{\partial^2 F}{\partial \kappa_i
\partial \kappa_j}(\lambda_a,p_a,q_a,\kappa_a)\right]=0.
$
\end{prop}
\begin{proof}
By Proposition \ref{regular} if $(\lambda_a,p_a,q_a)$ is a regular
point of $\mathbb E(L)$, the map (\ref{maprank})
has maximal rank $2m-1$. Also,
$\text{rank}\left[\frac{\partial^2 F}{\partial \kappa_i
\partial \kappa_j}(\lambda_a,p_a,q_a,\kappa_a)\right]$ is $2m-2$
which implies one of the columns of this matrix is linearly dependent on the others.
Assume this is the first column. Thus,
$\kappa \mapsto
\left(\frac{\partial F}{\partial \kappa_{[2m-1]\setminus
[1]}}(\lambda_a,p_a,q_a,\kappa), \ \det\left[\frac{\partial^2
F}{\partial \kappa_i
\partial \kappa_j}(\lambda_a,p_a,q_a,\kappa)\right]\right)$
has maximal rank. By implicit function theorem there is a
smooth map germ $\mathcal K:\mathbb R_e^{2m+1}\rightarrow \mathbb
R^{2m-1}$ at $(\lambda_a,a)$, s.t. $\kappa=\mathcal
K(\lambda,p,q)$ iff
$\frac{\partial
F}{\partial \kappa_{[2m-1]\setminus [1]}}(\lambda,p,q,\kappa)=0, \
\det\left[\frac{\partial^2 F}{\partial \kappa_i
\partial \kappa_j}(\lambda,p,q,\kappa)\right]=0.$
Then the germ of $\mathbb E(L)$ at$(\lambda_a,a)$ is
$
\mathbb E(L)=\left\{(\lambda,p,q):\frac{\partial F}{\partial
\kappa_{1}}(\lambda,p,q,\mathcal K(\lambda,p,q))=0\right\}.
$
The fiber of $\pi_r$ is tangent to $\mathbb E(L)$ at
$(\lambda_a,a)$ iff
\begin{equation}\label{df-d1}
\frac{\partial^2 F}{\partial \lambda
\partial \kappa_{1}}(\lambda_a,p_a,q_a,\kappa_a)+\sum_{j=1}^{2m-1}\frac{\partial^2 F}{\partial
\kappa_{j}
\partial \kappa_{1}}(\lambda_a,p_a,q_a,\kappa_a)\frac{\partial \mathcal K_j}{\partial \lambda}(\lambda_a,p_a,q_a)=0.
\end{equation}
Differentiating $\frac{\partial F}{\partial
\kappa_{[2m-1]\setminus [1]}}(\lambda,p,q,\mathcal
K(\lambda,p,q))=0$ w.r.t. $\lambda$ we obtain
\begin{equation}\label{df-di}
\frac{\partial^2 F}{\partial \lambda
\partial \kappa_{i}}(\lambda_a,p_a,q_a,\kappa_a)+\sum_{j=1}^{2m-1}\frac{\partial^2 F}{\partial
\kappa_{j}
\partial \kappa_{i}}(\lambda_a,p_a,q_a,\kappa_a)\frac{\partial \mathcal K_j}{\partial \lambda}(\lambda_a,p_a,q_a)=0.
\end{equation}
Thus (\ref{df-d1})-(\ref{df-di}) imply (\ref{rank}). But also (\ref{df-di}) and (\ref{rank}) imply (\ref{df-d1}).
\end{proof}
\begin{thm} \label{criminant}
The point $a=\lambda a^+ +(1-\lambda) a^-$ belongs to the
criminant $\Delta(L)$ of $GCS(L)$ iff there is a bitangent hyperplane to $L$ at
$a^+$ and $a^-$.
\end{thm}
\begin{proof}
If $(\lambda,a)\in\mathbb
E(L)$ is regular, by Propositions \ref{regular}-\ref{tangent}, $a^+, a^-$ are
$1$-parallel  and $a=(p,q)\in\Delta(L)$ iff
$(\lambda,a)$ satisfies (\ref{rank}). Thus
$\left[\frac{\partial^2
F}{\partial \kappa_i
\partial \kappa_j}\right] =$
$$
\frac{1}{2}\left[\begin{array}{ccccccc} \frac{\partial^2
S^+}{(\partial q_1^+)^2}-\frac{\partial^2 S^-}{(\partial
q_1^-)^2}& \frac{\partial^2 S^+}{\partial q_1^+ \partial
q_2^+}&\cdots &\frac{\partial^2 S^+}{\partial q_1^+ \partial
q_{m}^+}& -\frac{\partial^2 S^-}{\partial q_1^- \partial
p_2^-}&\cdots &-\frac{\partial^2 S^-}{\partial
q_1^- \partial p_{m}^-}\\
\frac{\partial^2 S^+}{\partial q_1^+ \partial q_2^+}&
\frac{\partial^2 S^+}{(
\partial q_2^+)^2}&\cdots &\frac{\partial^2 S^+}{\partial q_2^+
\partial q_{m}^+}& -1&\cdots &0\\
\vdots&\vdots&\ddots&\vdots&\vdots&\ddots&\vdots\\
\frac{\partial^2 S^+}{\partial q_1^+ \partial q_{m}^+}&
\frac{\partial^2 S^+}{
\partial q_2^+\partial q_{m}^+}&\cdots &\frac{\partial^2 S^+}{(\partial q_{m}^+
)^2}& 0&\cdots &-1\\
-\frac{\partial^2 S^-}{\partial q_1^- \partial
p_2^-}&-1&\cdots&0&-\frac{\partial^2 S^-}{(
\partial p_2^-)^2}&\cdots&\frac{\partial^2 S^-}{\partial p_2^-
\partial p_{m}^-}\\
\vdots&\vdots&\ddots&\vdots&\vdots&\ddots&\vdots\\
-\frac{\partial^2 S^-}{\partial q_1^- \partial
p_{m}^-}&0&\cdots&-1&\frac{\partial^2 S^-}{\partial p_2^-
\partial p_{m}^-}&\cdots&-\frac{\partial^2 S^-}{(
\partial p_{m}^-)^2}\\
\end{array}\right]
$$
and $
\ \frac{\partial^2 F}{\partial \lambda
\partial \beta_1}=p_1^+-p_1^--\sum_{j=1}^n q_j^+\frac{\partial^2 S^+}{\partial q_1^+ \partial
q_j^+}+q_1^-\frac{\partial^2 S^-}{(\partial q_1^-)^2}+\sum_{j=2}^n
p_j^-\frac{\partial^2 S^-}{\partial q_1^- \partial p_j^-},
$
$
\frac{\partial^2 F}{\partial \lambda
\partial \beta_i}=p_i^+-\sum_{i=1}^n q_j^+\frac{\partial^2 S^+}{\partial q_i^+ \partial
q_j^+}$,
$
\frac{\partial^2 F}{\partial \lambda
\partial \alpha_i}=q_i^-+q_1^-\frac{\partial^2 S^+}{\partial p_i^- \partial q_1^-
}+\sum_{j=2}^n p_j^-\frac{\partial^2 S^+}{\partial p_i^-
\partial p_j^-}$, for $ i=2,\cdots,m,$,
with $q^+=\frac{q+\beta}{2\lambda}$, $p^+=\frac{\partial
S^+}{\partial q^+}$  and
$q_1^-=\frac{q_1-\beta_1}{2(1-\lambda)}$, $p_{[m]\setminus
[2]}^-=\frac{p_{[m]\setminus [2]}-\alpha_{[m]\setminus
[2]}}{2(1-\lambda)}$ , $p_1^-=\frac{\partial S^-}{\partial
q_1^-}$, $q_{[m]\setminus [2]}^-=-\frac{\partial S^-}{\partial
p_{[m]\setminus [2]}^-}$.
Then, (\ref{rank}) is equivalent to
\begin{equation}\label{bitangent}
(a^+-a^-)\in T_{a^+}L^+ + T_{a^-}L^-,
\end{equation}
 since
$T_{a^+}L^+$ is spanned by  $\sum_{j=1}^m \frac{\partial^2
S^+}{\partial q_i^+\partial q_j^+}\frac{\partial}{\partial
p_j}+\frac{\partial}{\partial q_i}$ for $i=1,\cdots,m$ and
$T_{a^-}L^-$ is spanned by  $\frac{\partial^2
S^-}{(\partial q_1^-)^2}\frac{\partial}{\partial p_1}-\sum_{j=2}^m
\frac{\partial^2 S^-}{\partial q_1^-\partial
p_j^-}\frac{\partial}{\partial q_j}+\frac{\partial}{\partial q_1}$
and $\frac{\partial^2 S^-}{\partial p_i^-\partial
q_1^-}\frac{\partial}{\partial p_1}-\sum_{j=2}^m \frac{\partial^2
S^-}{\partial p_i^-\partial p_j^-}\frac{\partial}{\partial
q_j}+\frac{\partial}{\partial p_i}$ for $i=2,\cdots,m$.
If $a^+$, $a^-$ is $1$-parallel, (\ref{bitangent}) means  there is a bitangent hyperplane to $L^+$ at $a^+$ and to
$L^-$ at $a^-$. By continuity, a point in the closure of the set
of points satisfying (\ref{bitangent}) also satisfies this condition.
\end{proof}
\begin{cor} If, for some $\lambda$, $\lambda a^+ +(1-\lambda) a^- = a\in \Delta(L)\subset GCS(L)$, then the whole chord
$l(a^+,a^-)\subset GCS(L)$. Equivalently, if there is a
bitangent hyperplane to $L$ at $a^+$ and $a^-$, then $l(a^+,a^-)\subset GCS(L)$.
\end{cor}

Thus, we generalize the notion of
convexity of a curve on the plane.
\begin{defn}\label{weakconv}
A smooth closed Lagrangian submanifold $L$ of $(\mathbb R^{2m},\omega)$  is {\bf weakly convex}
if there is no bitangent hyperplane to $L$.
\end{defn}
\begin{cor}
If $L$ is a weakly convex closed Lagrangian submanifold of
$(\mathbb R^{2m},\omega)$  then the criminant $\Delta(L)$ of
$GCS(L)$ is empty.
\end{cor}

\section{Affine-Lagrangian stable singularities of the GCS}\label{section-sing-GCS}

We now define an equivalence relation to classify the singularities of $GCS(L)$. Due to the Lagrangian condition, we look for an equivalence of  generating families.
For the classification of $\mathbb E(\lambda)$ and $GCS(L)$, because $\lambda$ is no longer fixed it has become an extra parameter that unfolds the generating families $F$.
The naive approach is to consider the extended parameter space $\mathbb R\times\mathbb R^{2m}\ni (\lambda,x)$ for unfolding the generating families $f(\lambda,\kappa)=f_{\lambda}(\kappa)$ and classify their stable unfoldings in the usual way.
However, such a classification of $GCS(L)$ would not take into account the
projection $\pi: \mathbb R\times\mathbb R^{2m}\to\mathbb R^{2m}$
 in a proper way, because it  does not distinguish the affine time $\lambda\in\mathbb R$ from $x\in\mathbb R^{2m}$.

Now,  if $\mathcal A=(A,a)$ is an element of the affine symplectic group $iSp^{2m}_{\mathbb R}=Sp(2m,\mathbb R)\ltimes\mathbb R^{2m}$, with $A\in Sp(2m,\mathbb R)$, $a\in \mathbb R^{2m}$, then
\begin{equation}\label{A}
\mathcal A:(\mathbb R^{2m},\omega)\supset L\to L'\subset(\mathbb R^{2m},\omega) \ , \ x\mapsto \mathcal Ax=Ax+a \ .
\end{equation}
From this, we define the natural action
$$id_{T^*\mathbb R}\times\mathcal A\times\mathcal A:T^*\mathbb R\times\mathbb R^{2m}\times\mathbb R^{2m}\to T^*\mathbb R\times\mathbb R^{2m}\times\mathbb R^{2m} \ ,$$
$$(\lambda,\lambda^*,x^+,x^-)\mapsto(\lambda,\lambda^*,\mathcal Ax^+,\mathcal Ax^-) \ ,$$
which, via the chord transformation $\Phi_{\lambda}$,  induces an action
$$iSp^{2m}_{\mathbb R} \ \ni \ id_{T^*\mathbb R}\times\mathcal A_{\Phi} \ :  \  T^*\mathbb R\times T\mathbb R^{2m}\supset\mathcal L \ \to \  \mathcal L'\subset T^*\mathbb R\times T\mathbb R^{2m},$$
\begin{equation}\label{affine-action-TM}
id_{T^*\mathbb R}\times\mathcal A_{\Phi}: (\lambda,\lambda^*, x, \dot x)\mapsto (\lambda,\lambda^*, Ax+a, A\dot x + (2\lambda -1)a),
\end{equation}
that  commutes with projection $id_{T^*\mathbb R}\times pr: T^*\mathbb R\times T\mathbb R^{2m}\to T^*\mathbb R\times\mathbb R^{2m}$, that is,
defining the obvious action  $id_{\mathbb R}\times \mathcal A$ on $\mathbb R\times\mathbb R^{2m}$, we have
\begin{equation}\label{affine-invariance-of-pr}
(id_{\mathbb R}\times\mathcal A)\circ (id_{T^*\mathbb R}\times pr) = (id_{T^*\mathbb R}\times pr)\circ(id_{T^*\mathbb R}\times\mathcal A_{\Phi}).
\end{equation}
\begin{defn}\label{$(1,2m)$-Lagrangian-eq}
Germs of Lagrangian submanifolds $\mathcal L, \
\widetilde{\mathcal L}$ of the fiber bundle
$(T^{\ast}\mathbb R \times T\mathbb R^{2m}, d\lambda^*\wedge
d\lambda + \dot\omega)$ are {\bf (1,2m)-Lagrangian equivalent}
if there exists a symplectomorphism-germ $\Upsilon$ of
$T^{\ast}\mathbb R \times T\mathbb R^{2m}$ such that
$\Upsilon(\mathcal L)=\widetilde{\mathcal L}$ and the following
diagram commutes:
$$
\begin{array}{ccccccccc}

& & & Pr &  & \pi & \\
\mathcal L& \hookrightarrow &T^{\ast}\mathbb R\times T\mathbb
R^{2m} & \longrightarrow & \mathbb R\times \mathbb R^{2m} & \to &
\mathbb
R^{2m}  \\
& & & &  &  \\
& & \downarrow \Upsilon &  & \downarrow &  & \downarrow    \\
& & & Pr &  & \pi &  \\
\widetilde{\mathcal L}& \hookrightarrow &T^{\ast}\mathbb R\times
T\mathbb R^{2m} &  \longrightarrow & \mathbb R\times \mathbb
R^{2m} & \to & \mathbb
R^{2m}  \\
\end{array}
$$
The first two vertical diffeomorphism-germs (from right to left)
read:
$$ x\mapsto X(x) \ , \  (\lambda,x)\mapsto (\Lambda(\lambda,x),X(x)). $$
Moreover, germs $\mathcal L, \ \widetilde{\mathcal L}$ at
$(\frac{1}{2},a,\dot{a})$ are {\bf (1,2m)-Lagrangian equivalent
for $\lambda=\frac{1}{2}$ } if, in addition, for every $x\in
\mathbb R^{2m}$
\begin{equation}\label{L1/2}
\Lambda(\frac{1}{2},x)=\frac{1}{2}.
\end{equation}
\end{defn}
\begin{rem} \label{r1/2} Condition (\ref{L1/2}) is introduced for the
classification of the Wigner caustic $E_{{1}/{2}}(L)$ as a
part of $GCS(L)$. \end{rem}
\begin{defn}\label{GCS-$(1,2m)$-Lagrangian-eq}
$GCS(L)$ and $GCS(\widetilde{L})$ are {\bf (1,2m)-Lagrangian equivalent} if  $\mathcal L$ and
$\widetilde{\mathcal L}$ are (1,2m)-Lagrangian equivalent.
\end{defn}
\begin{rem}\label{GCS-affine-symplectic-eq} From (\ref{affine-invariance-of-pr}), it is clear that
 classification of $GCS(L)$ by $(1,2m)$-Lagrangian
equivalence of $\mathcal L$ is {\it affine symplectic invariant}.
\end{rem}
\begin{rem}\label{Lag-bif}
$(1,2m)$-Lagrangian equivalence of germs of Lagrangian
submanifolds of $(T^{\ast}\mathbb R
\times T\mathbb R^{2m}, d\lambda^*\wedge d\lambda + \dot\omega)$
is the equivalence of bifurcations of Lagrangian maps \cite{AGV}, that is, diagrams of the form:
$$
\begin{array}{ccccccccc}\label{diagram}
& & & & Pr &  & \pi & \\
D(\mathcal L):& \mathcal L& \hookrightarrow &T^{\ast}\mathbb
R\times T\mathbb R^{2m} & \longrightarrow & \mathbb R\times
\mathbb R^{2m} & \to &
\mathbb R^{2m}\\
\end{array}
$$
\end{rem}
\begin{defn}
$\mathcal L$ is {\bf (1,2m)-Lagrangian
stable}  if the diagram of maps $D(\mathcal L)$ is stable, i.e.
every $\widetilde{\mathcal L}$ with nearby
diagram $D(\widetilde{\mathcal L})$ is $(1,2m)$-Lagrangian
equivalent  to $\mathcal L$. $GCS(L)$ is {\bf (1,2m)-Lagrangian
stable} if $\mathcal L$ is (1,2m)-Lagrangian stable. In view of Remark \ref{GCS-affine-symplectic-eq}, we also use
the term {\bf affine-Lagrangian stability} for $(1,2m)$-Lagrangian
stability.
\end{defn}
\begin{defn}\label{def-gf}
The function-germs $F, \widetilde{F}:\mathbb R\times \mathbb R^{2m}
\times \mathbb R^k\rightarrow \mathbb R$ are {\bf (1,2m)-$\mathcal
R^+$-equivalent} if there exists a diffeomorphism-germ
$$(\lambda,x,\kappa)\mapsto(\Lambda(\lambda,x),X(x),K(\lambda,x,\kappa))$$ and a
smooth function-germ $g:\mathbb R\times \mathbb R^{2m}\rightarrow
\mathbb R$ such that
$$
\widetilde{F}(\lambda,x,\kappa)=F(\Lambda(\lambda,x),X(x),K(\lambda,x,\kappa))+g(\lambda,x).
$$
Germs $F$ and $\widetilde F$ with the common $(\lambda,x)$-space
$\mathbb R\times \mathbb R^{2m}$ of parameters, but in general,
with their spaces of arguments of different dimensions are {\bf
stably (1,2m)-$\mathcal R^+$-equivalent} if there are
nondegenerate quadratic forms $Q$ in new arguments $\xi$ and
$\widetilde Q$ in new arguments $\tilde \xi$ such that $F+Q$ and
$\widetilde F + \widetilde Q$ are $(1,2m)$-$\mathcal
R^+$-equivalent. The germ $F$ at $(\frac{1}{2},a,\kappa_a)$ and
the germ $\widetilde F$ at $(\frac{1}{2},a,\tilde \kappa_a)$ are
 (stably) (1,2m)-$\mathcal R^+$-equivalent for
$\lambda=\frac{1}{2}$ if, in addition, for every $x\in \mathbb
R^m$ condition (\ref{L1/2}) is satisfied.
\end{defn}
\begin{rem}
$(1,2m)$-$\mathcal R^+$-equivalence is a special case of Wassermann's
$(1,2m)$-equivalence  \cite{Was}. For relations between $(r,s)$-classification
of families of functions \cite{Was}, classification of
bifurcations of caustics \cite{A-b}, \cite{Z} and
classification of bifurcations of Lagrangian maps, see \cite{AGV}.
\end{rem}

 We have the following result, whose proof is a minor
modification for $(1,2m)$-Lagrangian equivalence of the proof of
Theorem \ref{LagstabR+} in \cite{AGV}.
\begin{prop}\label{1-2m-R-Lag}
Germs of Lagrangian submanifolds $\mathcal L, \
\widetilde{\mathcal L}$ of
$(T^{\ast}\mathbb R \times T\mathbb R^{2m}, d\lambda^*\wedge
d\lambda + \dot\omega)$ are $(1,2m)$-Lagrangian equivalent iff the germs of generating families $F$ and
$\widetilde{F}$ are stably $(1,2m)$-$\mathcal R^+$-equivalent.
\end{prop}
\begin{defn}\label{GCS-affine-Lag-stability}
The function-germ $F$ at $z$ is {\bf (1,2m)-$\mathcal R^+$-stable}
if  $\forall$ neighborhood $U\ni z$ in $\mathbb R\times \mathbb
R^{2m} \times \mathbb R^k$ and  representative function $F'$ of
the germ $F$ on $U$, there is a neighborhood $V$ of
$F'$ in $C^{\infty}(U,\mathbb R)$ (with weak
$C^{\infty}$-topology) such that for any function $G'\in V$ there
is a point $z'\in U$ such that the germ of $G'$ at $z'$ is
$(1,2m)$-$\mathcal R^+$-equivalent to $F$.
\end{defn}
\begin{rem}  $\mathcal L$ and $GCS(L)$ are $(1,2m)$-Lagrangian stable if and only if  the germ of generating family
 $F$  (of $\mathcal L$) is $(1,2m)$-$\mathcal R^+$-stable.
\end{rem}
The following theorems
show that the only affine-Lagrangian stable singularities of GCS
are singularities of the criminant, the smooth part of the Wigner
caustic and their ``tangent'' union.
\begin{thm} \label{class-criminant}
Let $\lambda_a\ne \frac{1}{2}$. If $F$ is the germ at $(\lambda_a,a,\kappa_a)$ of a $(1,2m)$-$\mathcal
R^+$-stable unfolding of $f\in \mathfrak m^2$ then $F$ is stably
$(1,2m)$-$\mathcal R^+$-equivalent to the germ of the trivial
unfolding or to one of the
following germs at $(0,0,0)$ of unfoldings  of $f(t)=t^3$
\begin{equation}\label{A2Ak}
A_2^{A_k^{\pm}} \ : \
F(\lambda,x,t)=t^3+t\left(\sum_{i=1}^{k}x_i\lambda^{i-1}\pm\lambda^{k+1}\right),
\end{equation}
for $k=0,1,2,\cdots,2m$ (the notation $A_2^{A_k^{\pm}}$ is taken from \cite{Go}).
\end{thm}
\begin{proof}
If $f$ has $A_1$ singularity then it is obvious that $F$ is stably
$(1,2m)$-$\mathcal R^+$-equivalent to the trivial unfolding. Now
we assume that $f$ has $A_2$ singularity. Since $F$ is stable then
$F$ is stable $(1,2m)$-$\mathcal R^+$-equivalent to
$F(\lambda,x,t)=t^3+t g(\lambda,x)$, where $g$ is a smooth
function-germ vanishing at $0$. If $g$ is a versal unfolding of
the function-germ $\lambda\mapsto g(\lambda,0)$ with $A_k$
singularity we can reduce $F$ to the form (\ref{A2Ak}) by a
diffeomorphism-germ of the form $(\lambda,x,t)\mapsto
(\Lambda(\lambda,x),X(x),t)$.
The following lemma shows that these are the only $(1,2m)$-$\mathcal R^+$-stable
unfoldings.
\end{proof}
\begin{lem}
Unfoldings of $A_3^{\pm}$ singularity are not $(1,2m)$-$\mathcal
R^+$-stable.
\end{lem}
\begin{proof}  If $f$ has $A_3$ singularity then
$F$ is stable $(1,2m)$-$\mathcal R^+$-equivalent to
$F(\lambda,x,t)=\pm t^4+t^2 g_2(\lambda,x)+tg_1(\lambda,x)$, where
$g_1, g_2$ are smooth function-germs vanishing at $0$. Now we use
the standard arguments of the singularity theory that stability
implies infinitesimal stability. In the case of $(1,2m)$-$\mathcal
R^+$-equivalence, the infinitesimal stability implies
\begin{equation}\label{inf-st-R+}
\mathcal E_2=\mathcal E_2\left<\frac{\partial F}{\partial
t}|_{\mathbb R^2}\right>+\mathcal E_1\left<1,\frac{\partial
F}{\partial \lambda}|_{\mathbb R^2}\right>+\mathbb
R\left<\frac{\partial F}{\partial x_1}|_{\mathbb
R^2},\cdots,\frac{\partial F}{\partial x_{2m}}|_{\mathbb
R^2}\right>+\mathfrak m_2^{2m+4},
\end{equation}
where $\mathbb R^2$ denotes the $t,\lambda$-plane $\{x=0\}$,
$\mathcal E_2$ is the ring of smooth function-germs in $\lambda$
and $t$, $\mathfrak m_2$ is the maximal ideal in $\mathcal E_2$
and $\mathcal E_1$ is the ring of smooth function-germs in
$\lambda$. Now we use the method from \cite{Was}.
Let $V=\mathcal E_2\left/\left(\mathcal E_2\left<\frac{\partial
F}{\partial t}|_{\mathbb R^2}\right>+\mathfrak
m_2^{2m+4}\right)\right.$ and let $\pi:\mathcal E_2\rightarrow V$. We have $\pi(t^3)=\pi(\mp 1/2 tg_2|_{\mathbb
R^2}\mp 1/4 g_1|_{\mathbb R^2})$ in $V$. Thus  elements $\pi(t^i
\lambda^j)$ for $i=0,1,2$ and $j<2m+4-i$ form a basis of $V$ over
$\mathbb R$. Thus, $\dim _{\mathbb R} V=6m+9$.
   Moreover, $\frac{\partial F}{\partial
\lambda}|_{\mathbb R^2}=t\left(t\frac{\partial g_2}{\partial
\lambda}|_{\mathbb R^2}+\frac{\partial g_1}{\partial
\lambda}|_{\mathbb R^2}\right)$. Then $ \dim_{\mathbb
R}\pi\left(\mathcal E_1\left<1,\frac{\partial F}{\partial
\lambda}|_{\mathbb R^2}\right>\right)\le 4m+7$,  $
\dim_{\mathbb R}\pi\left(\mathbb R\left<\frac{\partial F}{\partial
x_1}|_{\mathbb R^2},...,\frac{\partial F}{\partial
x_{2m}}|_{\mathbb R^2}\right>\right)\le 2m.
$
So, (\ref{inf-st-R+}) implies $\dim_{\mathbb R}V\le
6m+7<6m+9$, which is impossible. Thus $F$ is not
$(1,2m)$-$\mathcal R^+$-stable, $A_3$ has no
$(1,2m)$-$\mathcal R^+$-stable unfoldings.
\end{proof}
For $E_{1/2}(L)\subset GCS(L)$,  we consider the germ of
$F$ at $(1/2,a,\kappa_a)$.
\begin{thm} \label{class-Wigner} If $F$ is the germ at $(\frac{1}{2},a,\kappa_a)$ of a $(1,2m)$-$\mathcal
R^+$-stable unfolding of $f\in \mathfrak m^2$ then $F$ is stably
$(1,2m)$-$\mathcal R^+$-equivalent ($\lambda=1/2$) to the germ
of the trivial unfolding  or to one
of the following germs at $(\frac{1}{2},0,0)$ of unfoldings  of
$f(t)=t^3$
\begin{equation}\label{A2Ak1/2}
A_2^{B_k^{\pm}} \ : \
F(\lambda,x,t)=t^3+t\left(\sum_{i=0}^{k-1}x_{i+1}\left(\lambda-\frac{1}{2}\right)^{i}\pm\left(\lambda-\frac{1}{2}\right)^{k}\right),
\end{equation}
for $k=1,2,\cdots,2m$ (the notation $A_2^{B_k^{\pm}}$ is taken from \cite{Go}).
\end{thm}
\begin{proof}
If $f$ has $A_1$ singularity then  $F$ is stably
$(1,2m)$-$\mathcal R^+$-equivalent to the trivial unfolding. If $f$ has $A_2$ singularity, since $F$ is stable then
$F$ is stably $(1,2m)$-$\mathcal R^+$-equivalent to
$F(\lambda,x,t)=t^3+t g(\lambda,x)$, where $g$ is a smooth
function-germ vanishing at $(1/2,0)$. If $g$ is a versal unfolding
of the function-germ $\lambda\mapsto g(\lambda,0)$ with
$B_k^{\pm}$ singularity on a manifold ($\lambda$-space) with  the
boundary ($\lambda=\frac{1}{2}$) (see \cite{A-b}) then we can
reduce $F$ to the form (\ref{A2Ak1/2}) by a diffeomorphism-germ of
the form $(\lambda,x,t)\mapsto
(1/2+(\lambda-1/2)\Lambda(\lambda,x),X(x),t)$.
\end{proof}
\begin{thm} \label{sing-criminant}
If $F$ (generating  $\mathcal L$) has
$A_2^{A_k^\pm}$ singularity, for $k=0,1,\cdots, 2m$, then $\mathbb
E(L)$ is a germ of a smooth hypersurface in $\mathbb R\times
\mathbb R^{2m}$.

If $F$ has $A_2^{A_0}$ singularity at $(\lambda_a,a,\kappa_a)$
then $\mathbb E(L)$ is transversal at $(\lambda_a,a)$ to the
fibers of projection $\pi$.

If $F$ has $A_2^{A_k^{\pm}}$ singularity for $k\ge 1$ at
$(\lambda_a,a,\kappa_a)$ then $\mathbb E(L)$ is $k$-tangent at
$(\lambda_a,a)$ to the fibers of $\pi$, $a$ belongs to
the criminant $\Delta(L)$ of $GSC(L)$ and the germ of $\Delta(L)$
at $a$ is the caustic of $A_k^{\pm}$ singularity.
\end{thm}
\begin{proof}
By Proposition \ref{wl} and the normal form of $F$ for
$A_2^{A_k^{\pm}}$ singularity we obtain
$
\mathbb E(L)=\{(\lambda,x)\in \mathbb R\times \mathbb
R^{2m}:\sum_{i=1}^{k}x_i\lambda^{i-1}\pm\lambda^{k+1}=0 \}$.
It is easy to see that $\mathbb E(L)$ is the germ at $(0,0)$ of a
smooth hypersurface and $\mathbb E(L)$ is transversal at $(0,0)$
to $\{\lambda=0\}$ for $k=0$ and $\mathbb E(L)$ is $k$-tangent to
$\{\lambda=0\}$ at $(0,0)$ for $k=1,2,\cdots,2m$. The germ of $\Delta(L)$ at $0$  is
$$
\{x\in \mathbb R^{2m}:\exists \lambda \
\sum_{i=1}^{k}x_i\lambda^{i-1}\pm\lambda^{k+1}=0, \
\sum_{i=2}^{k}(i-1)x_i\lambda^{i-2}\pm(k+1)\lambda^{k}=0 \}.
$$
So $\Delta(L)$ is a caustic of $A_k^\pm$ singularity.
\end{proof}
\begin{thm} \label{sing-Wigner}
If the germ at $(\frac{1}{2},a,\kappa_a)$ of
$F$  has $A_2^{B_k^\pm}$ singularity
($k=1,\cdots, 2m$), then $\mathbb E(L)$ is a germ of smooth
hypersurface in $\mathbb R\times \mathbb R^{2m}$.

If $F$ has $A_2^{B_1}$ singularity at
$(\frac{1}{2},a,\kappa_a)$ then $\mathbb E(L)$ is transversal at
$(\frac{1}{2},a)$ to the fibers of projection $\pi$. The germ of
$GCS(L)$ at $a$ is the germ of a smooth hypersurface of $\mathbb
R^{2m}$ - the Wigner caustic $E_{{1}/{2}}(L)$.

If $F$ has $A_2^{B_k^{\pm}}$ singularity for $k\ge 2$ at
$(\frac{1}{2},a,\kappa_a)$ then $\mathbb E(L)$ is $k$-tangent at
$(1/2,a,t)$ to the fibers of $\pi$. The germ of
$GCS(L)$ at $a$ consists of two tangent components: the germ of a
smooth hypersurface - $E_{{1}/{2}}(L)$ - and
the germ of the caustic of $B_k^{\pm}$ singularity -
$\Delta(L)$.
\end{thm}
\begin{proof}
By Proposition \ref{wl} and the normal form of $F$ for
$A_2^{B_k^{\pm}}$ singularity we get
$\mathbb E(L)=\{(\lambda,x)\in \mathbb R\times \mathbb
R^{2m}:\sum_{i=0}^{k-1}x_{i+1}(\lambda-1/2)^{i}\pm(\lambda-1/2)^{k}=0
\}.$ $E_{1/2}(L)=\{x\in \mathbb R^{2m}:x_1=0 \}$
is a germ of smooth hypersurface.
Thus $\mathbb E(L)$ is the germ at $(1/2,0)$ of
a smooth hypersurface and $\mathbb E(L)$ is transversal at
$(1/2,0)$ to $\{\lambda=1/2\}$ for $k=1$. For $k=2,\cdots,2m$,
$\mathbb E(L)$ is $k$-tangent to $\{\lambda=1/2\}$ at $(1/2,0)$.  The germ of  $\Delta(L)$ at $0$ is
$$
\{x\in \mathbb R^{2m}:\exists \tau \
\sum_{i=0}^{k-1}x_{i+1}\tau^{i}\pm\tau^{k}=0, \
\sum_{i=1}^{k-1}ix_{i+1}\tau^{i-1}\pm k\tau^{k-1}=0 \}.
$$
So $\Delta(L)$ is a caustic of $B_k^\pm$  and
$E_{1/2}(L)$ is tangent to $\Delta(L)$ at $0$.
\end{proof}
\begin{rem}\label{GCS-realization}
Not all $(1,2m)$-$\mathcal R^+$-stable singularities  can be
realizable as singularities of generating families $F$ for
$\mathcal L$ which are of the special form given in Theorem \ref{gen-fam}. In the next section, in Theorem \ref{real-lag-curve},
we prove that the $A_2^{A_2}$ singularity is not realizable for
Lagrangian curves.
\end{rem}

\section{Classifications of the GCS of Lagrangian curves}\label{css-class-Lag-curves}

We now classify the singularities of the global centre symmetry set of a Lagrangian curve
$L\subset(\mathbb R^2,\omega)$.
To set the stage, we first state the results for the $GCS$ of a curve on the affine plane $\mathbb R^2$,
when no symplectic structure is considered.
These results, obtained in  \cite{Ber},
\cite{Jan} and \cite{GZ1}-\cite{GZ2} by various methods,  are summarized in Theorem \ref{curvestability} and can also be proved using the affine-invariant method
of chord equivalence, the analogous of
$(1,2m)$-Lagrangian equivalence when no symplectic structure is
considered.
Theorem \ref{globalresults} presents global results for the $GCS$ of a convex curve,
some of which have not been stated before.

\begin{thm}[\cite{Ber}, \cite{Jan}, \cite{GZ1}-\cite{GZ2}]\label{curvestability}
Affine stable GCS of a
smooth convex closed curve $M\subset\mathbb R^2$ (no symplectic structure)
consists of: \\
i) The CSS, a smooth curve with (possible) self
intersections and cusp singularities, ii) the Wigner
caustic, a smooth curve with (possible) self intersections and
cusp singularities lying on the smooth part of the CSS,
and iii) the middle axes, which are smooth half-lines starting at
the the cusp points of the CSS.
\end{thm}

\begin{thm}\label{globalresults}
Let $M$ be a generic smooth convex closed curve in $\mathbb{R}^2$.
The number of cusps of the Wigner caustic of $M$ is odd and not
smaller than 3. The number of cusps of the CSS of $M$ is odd and
not smaller than 3. The number of cusps of the Wigner caustic of
$M$ is not greater than the number of cusps of the CSS of $M$.
\end{thm}

\begin{center}
\includegraphics{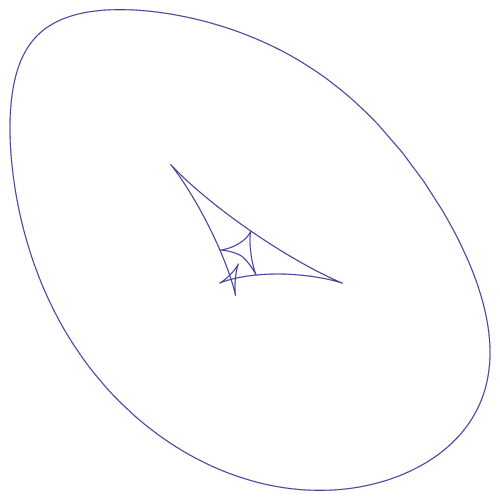}

{\small {\bf Figure 1.} GCS of an oval in the
plane: CSS with 5 cusps and Wigner caustic with 3 cusps
(the middle axes are not shown here).}
\end{center}

\begin{center}
\includegraphics{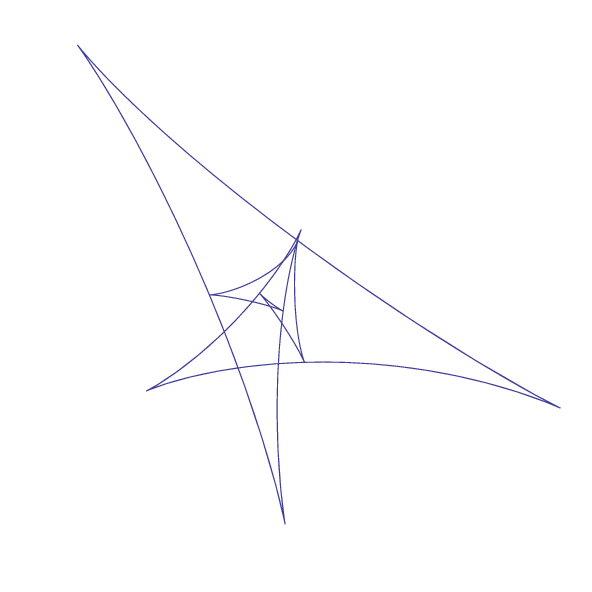}

{\small {\bf Figure 2.} Both the CSS and the Wigner caustic with
five cusps.}
\end{center}


\begin{proof}
The statement on the number of cusps of Wigner caustics,
was first proven by Berry \cite{Ber} and the statement on
the number of cusps of CSS, was first proven by Giblin and Holtom
\cite{GH}. The last inequality follows immediately from the
characterization in \cite{GH} of cusps of $E_{1/2}(M)$ by the
curvature ratio being $1$ and cusps of CSS of $M$ by the
derivative of the curvature ratio being $0$, using Rolle's
theorem.
\end{proof}
Figures of $GCS(M)$ where the number of cusps of the CSS and of
the Wigner caustic are equal to three and neither curve is self
intersecting can be found in \cite{GH}. We picture  a case
when the number of cusps of the Wigner caustic is three and the
CSS is self intersecting and the number of its cusps is five, and
another case when both the Wigner caustic and the CSS are self
intersecting and both have five cusps.

\subsection{Affine-Lagrangian  classification of the $GCS$ of Lagrangian curves}
Let  $L$ be a smooth closed curve on $(\mathbb R^2,\omega=dp\wedge dq)$. Using the
$(1,2)$-Lagrangian equivalence introduced in
Definition \ref{GCS-$(1,2m)$-Lagrangian-eq}, we classify the
singularities of $GCS(L)$.
In what follows,
$a^+=(p_a^+,q_a^+), a^-=(p_a^-,q_a^-)$ denote a parallel pair
on $L$ and $a_{\lambda}=\lambda a^++(1-\lambda)a^-$,
$\dot{q}_\lambda=\lambda q_a^+-(1-\lambda)q_a^-$. Let $S^{\pm}$ be
germs of generating functions of $L$ at $a^\pm$ satisfying the
conditions in Proposition \ref{s}. The germ of generating
family of $\mathcal L$ and the big wave front set are given by
$$
F(\lambda,p,q,t)=2\lambda^2S^+(\frac{q+t}{2\lambda})-2(1-\lambda)^2S^-(\frac{q-t}{2(1-\lambda)})-pt.
$$
$$
\mathbb E(L)=\left\{(\lambda,p,q)\in \mathbb R\times \mathbb R^2:
\exists t \ \frac{\partial F}{\partial t}(\lambda,p,q,t)=
\frac{\partial^2 F}{\partial t^2}(\lambda,p,q,t)=0 \right\}.
$$

The following propositions present  geometrical descriptions of
positions of $\mathbb E(L)$ with respect to  $\pi$ in terms of functions $F$,  $S^{+}$ and $S^{-}$.

\begin{prop} \label{curve-reg}
The following conditions are equivalent
\begin{itemize}
\item[(i)] $(\lambda,a_{\lambda})$ belongs the regular part of
$\mathbb E(L)$,

\item[(ii)] $\exists t \ \frac{\partial^3 F}{\partial
t^3}(\lambda,a_{\lambda},t)\ne 0,  \frac{\partial F}{\partial
t}(\lambda,a_{\lambda},t)=\frac{\partial^2 F}{\partial
t^2}(\lambda,a_{\lambda},t)=0$,

\item[(iii)] $\frac{1}{\lambda}\frac{\partial^3 S^+}{\partial
(q^+)^3}(q_a^+)+\frac{1}{1-\lambda}\frac{\partial^3 S^-}{\partial
(q^-)^3}(q_a^-)\ne 0$,

\item[(iv)]
$\frac{1}{\lambda}\kappa(a^+)+\frac{1}{1-\lambda}\kappa(a^-)\ne
0$, where $\kappa(x)$ is the curvature of $L$ at $x$.
\end{itemize}
\end{prop}
\begin{proof} Equivalence of (i) and (ii) follows from the
definition of the regular part of $\mathbb E(L)$. Equivalence of
(ii) and (iii) is obtained by direct calculations. (iv) is obvious
since $\kappa(a^{\pm})=\frac{\partial^3 S^\pm}{\partial
(q^\pm)^3}(q_a^\pm)$.
\end{proof}
\begin{prop}\label{curve-tangent}
The following conditions are equivalent
\begin{itemize}
\item[(v)] the regular part of $\mathbb E(L)$ is tangent to the
fiber of $\pi$ at $(\lambda,a_\lambda)$,

\item[(vi)] $\exists t$ (ii) is satisfied and $ \frac{\partial^2
F}{\partial \lambda
\partial t}(\lambda,a_{\lambda},t)=0$.

\item[(vii)] (iii) is satisfied and $p_a^+=\frac{\partial
S^+}{\partial q^+}(q_a^+)=\frac{\partial S^-}{\partial
q^-}(q_a^-)=p_a^-$.

\item[(viii)] (iv) is satisfied and $l(a^+,a^-)$ is bitangent to
$a^+,a^-$ to $L$.
\end{itemize}
\end{prop}
\begin{proof}
All statements follow from Proposition \ref{tangent} and Theorem
\ref{criminant}.
\end{proof}
\begin{prop}\label{1-tangent}
The following conditions are equivalent
\begin{itemize}
\item[(ix)] the regular part of $\mathbb E(L)$ is $1$-tangent to
the fiber of $\pi$ at $(\lambda,a_\lambda)$,

 \item[(x)] $\exists t$
(vi) is satisfied and
\begin{equation}\label{delta}
\left(\frac{\partial^3
F}{\partial \lambda
\partial t^2}(\lambda,a_\lambda,t)\right)^2-\frac{\partial^3 F}{
\partial t^3}(\lambda,a_\lambda,t)\frac{\partial^3 F}{\partial
\lambda^2 \partial t}(\lambda,a_\lambda,t)\ne 0.
\end{equation}

\item[(xi)] (vii) is satisfied and $ \frac{\partial^3
S^+}{\partial (q^+)^3}(q_a^+)\frac{\partial^3 S^-}{\partial
(q^-)^3}(q_a^-)\ne 0$.

\item[(xii)] (iv) is satisfied and $l(a^+,a^-)$ is $1$-tangent to
$L$ at $a^+$ and $a^-$
\end{itemize}
\end{prop}

\begin{proof} $(\lambda,a_\lambda)\in \mathbb E(L)$ is regular. By
Proposition \ref{curve-reg}, $\frac{\partial^3
F}{\partial t^3}(\lambda,a_\lambda,t)\ne 0$. Thus, $\exists$ smooth function-germ $T$ on $\mathbb R^3$ s.t.
$\frac{\partial^2 F}{\partial t^2}(\lambda,p,q,t)=0$ iff
$t=T(\lambda,p,q)$. Then $\mathbb E(L)=\left\{(\lambda,p,q):
\frac{\partial F}{\partial t}(\lambda,p,q,T(\lambda,p,q ))=0
\right\}$. Then
\begin{equation}\label{1-diff}
\frac{\partial}{\partial \lambda}\left(\frac{\partial F}{\partial
t}(\lambda,p,q,T(\lambda,p,q
))\right)\left|_{(\lambda,a_{\lambda})}\right.=0
\end{equation}
\begin{equation}\label{2-diff}
\frac{\partial^2}{\partial \lambda^2}\left(\frac{\partial
F}{\partial t}(\lambda,p,q,T(\lambda,p,q
))\right)\left|_{(\lambda,a_{\lambda})}\right.\ne 0
\end{equation}
are equivalent to (ix). Using the formula
\begin{equation}\label{ift-diff}
\frac{\partial T}{\partial
\lambda}(\lambda,p,q)=-\left(\frac{\partial^2 F}{
\partial t^3}(\lambda,p,q,T(\lambda,p,q)\right)^{-1}\frac{\partial^2 F}{\partial \lambda \partial
t^2}(\lambda,p,q,T(\lambda,p,q))
\end{equation}
we see that (\ref{1-diff})-(\ref{2-diff}) are
equivalent to (x). Equivalence of (x) and (xi) is obtained by
direct calculation and the last equivalence is obvious.
\end{proof}

\begin{prop}\label{2-tangent}
The following conditions are equivalent
\begin{itemize}
\item[(xiii)] the regular part of $\mathbb E(L)$ is $2$-tangent to
the fiber of $\pi$ at $(\lambda,a_\lambda)$,
 \item[(xiv)] $\exists t$
(vi) is satisfied, (\ref{delta}) is not satisfied and \\
$ \Big{\{}\frac{\partial^4 F}{\partial \lambda^3 \partial
t}\left(\frac{\partial^3 F}{
\partial t^3}\right)^3-3\frac{\partial^4 F}{\partial \lambda^2 \partial
t^2}\left(\frac{\partial^3 F}{
\partial t^3}\right)^2\frac{\partial^3 F}{
\partial \lambda \partial t^2}$\\
\noindent $+3\frac{\partial^4 F}{\partial \lambda \partial
t^3}\frac{\partial^3 F}{
\partial t^3}\left(\frac{\partial^3 F}{
\partial \lambda \partial t^2}\right)^2-\frac{\partial^4 F}{\partial
t^4}\left(\frac{\partial^3 F}{\partial \lambda
\partial t^2})^3 \right\} (\lambda,a_\lambda,t)\neq 0$
\item[(xv)] (vii) is satisfied and $\left(\frac{\partial^3 S^+}{\partial
(q^+)^3}(q_a^+)= 0 \  \wedge \  \frac{\partial^4 S^+}{\partial
(q^+)^4}(q_a^+)\ne 0\right)$ or \\ $\left(\frac{\partial^3 S^-}{\partial
(q^-)^3}(q_a^-)= 0 \ \wedge \ \frac{\partial^4 S^-}{\partial
(q^-)^4}(q_a^-)\ne 0\right)$
\item[(xvi)] (iv) is satisfied and $l(a^+,a^-)$ is $1$-tangent to
$L$ at one of points $a^+, a^-$ and $2$-tangent to $L$ at the
other.
\end{itemize}
\end{prop}
\begin{proof}
 (xiii) means that (\ref{1-diff}) is satisfied,
(\ref{2-diff}) is not satisfied and
$\frac{\partial^3}{\partial \lambda^3}\left(\frac{\partial
F}{\partial t}(\lambda,p,q,T(\lambda,p,q
))\right)\left|_{(\lambda,a_{\lambda})}\right.\ne 0$.
Using (\ref{ift-diff}), we see  that these conditions
are equivalent to (xiv). By direct calculation we see that
(xiv) $\iff$ (xv). Finally, (xvi) is the geometric
description of (xv).
\end{proof}
\begin{thm}\label{real-lag-curve}
Let $\frac{1}{\lambda}\frac{\partial^3 S^+}{\partial
(q^+)^3}(q_a^+)+\frac{1}{1-\lambda}\frac{\partial^3 S^-}{\partial
(q^-)^3}(q_a^-)\ne 0$ (for  (\ref{a01/2})-(\ref{a11/2}) below,
$\lambda=1/2$). Let $l(a^+,a^-)$ denote the chord passing through $(a^+,a^-)$.
\begin{enumerate}
\item \label{a01/2}If  $l(a^+,a^-)$ is not bitangent to
$L$ at $a^+,a^-$ then the germ of $F$ at
$(1/2,a_{1/2},\dot{q}_{1/2})$ has $A_2^{B_1}$ singularity and the
germ of $GCS$ at $a_{1/2}$ is a smooth curve (the smooth part of
the Wigner caustic).
\item  \label{a11/2} If $l(a^+,a^-)$ is $1$-tangent to
$L$ at $a^+$ and at $a^-$ then the germ of $F$ at
$(1/2,a_{1/2},\dot{q}_{1/2})$ has $A_2^{B_2}$ singularity and the
germ of $GCS$ at $a_{1/2}$ is a union of two $1$-tangent smooth
curves (the smooth part of the Wigner caustic and the smooth part
of the criminant).
\item \label{a1} If  $l(a^+,a^-)$ is $1$-tangent to $L$
at $a^+$ and at $a^-$ then the germ of $F$ at
$(\lambda,a_\lambda,\dot{q}_\lambda)$ for $\lambda\ne 1/2$ has
$A_2^{A_1}$ singularity and the germ of $GCS$ at $a_{\lambda}$ is
a smooth curve (the smooth part of the criminant).
\item \label{a2} If  $l(a^+,a^-)$ is $1$-tangent to $L$
at one of the points $a^+, a^-$ and $2$-tangent at the
other,  the germ of $F$ at
$(\lambda,a_\lambda,\dot{q}_\lambda)$ for $\lambda \ne 1/2$ is not
$(1,2)$-$\mathcal R^+$-stable.  $A_2^{A_2}$ is not realizable as stable
singularity of the $GCS$ of a Lagrangian curve.
\end{enumerate}
\end{thm}
\begin{proof}
By Proposition \ref{curve-reg}, if
$\frac{1}{\lambda}\frac{\partial^3 S^+}{\partial
(q^+)^3}(q_a^+)+\frac{1}{1-\lambda}\frac{\partial^3 S^-}{\partial
(q^-)^3}(q_a^-)\ne 0$
then the germ of
$F$ is a unfolding of $A_2$ singularity. Therefore we can reduce $F$
to the  form $F'(\lambda,p,q,t)=t^3+g(\lambda,p,q)t$,
where $g$ is a smooth function-germ vanishing at $(\lambda_a,0)$
(for $\lambda_a=0$ or $\lambda_a=1/2$).
By Proposition \ref{curve-tangent}, if $l(a^+,a^-)$ is
not bitangent to $L$ at $a^+, a^-$ then $\frac{\partial
F'}{\partial t
\partial \lambda}(1/2,0,0) \ne 0$ and this implies
$\frac{\partial g}{\partial \lambda}(1/2,0)\ne 0$. By Theorems
\ref{class-Wigner} and \ref{sing-Wigner} we obtain (\ref{a01/2}).
If the chord $l(a^+,a^-)$ is tangent to $L$ at $a^+, a^-$ then by
Proposition \ref{curve-tangent} we get that $p^+_a=p^-_a$ and
$\frac{\partial F'}{\partial t
\partial \lambda}(\lambda_a,0,0) = 0$ and this implies
$\frac{\partial g}{\partial \lambda}(\lambda_a,0)=0$. But
$dg|_{(\lambda_a,0)}\ne 0$ since $\frac{\partial F}{\partial t
\partial p}(\lambda_a,a,\dot{q}_a) \ne 0$.
By Proposition \ref{1-tangent} if  $l(a^+,a^-)$ is $1$-tangent to
$L$ at $a^+, a^-$ then
$\left(\frac{\partial^3 F'}{\partial \lambda \partial
t^2}(\lambda_a,0,0)\right)^2-\frac{\partial^3 F'}{
\partial t^3}(\lambda_a,0,0)\frac{\partial^3 F'}{\partial
\lambda^2 \partial t}(\lambda_a,0,0)\ne 0$.
But this implies $\frac{\partial^2 g}{\partial
\lambda^2}(\lambda_a,0,)\ne 0$. Thus if $\lambda_a=1/2$ by
Theorems \ref{class-Wigner} and \ref{sing-Wigner} we obtain
(\ref{a11/2}) and otherwise by Theorems \ref{class-criminant} and
\ref{sing-criminant} we obtain (\ref{a1}).
Finally, assume that  $l(a^+,a^-)$ is $1$-tangent
to $L$ at $a^+$ and $2$-tangent at $a^-$. By Proposition
\ref{2-tangent} we get  $\frac{\partial^2 g}{\partial
\lambda^2}(\lambda_a,0,)=0$ and
$\Big\{\frac{\partial^4 F}{\partial
\lambda^3 \partial t}\left(\frac{\partial^3 F}{
\partial t^3}\right)^3-3\frac{\partial^4 F}{\partial \lambda^2 \partial
t^2}\left(\frac{\partial^3 F}{
\partial t^3}\right)^2\frac{\partial^3 F}{
\partial \lambda \partial t^2}+3\frac{\partial^4 F}{\partial \lambda \partial
t^3}\frac{\partial^3 F}{
\partial t^3}\left(\frac{\partial^3 F}{
\partial \lambda \partial t^2}\right)^2-\frac{\partial^4 F}{\partial
t^4}\left(\frac{\partial^3 F}{\partial \lambda
\partial t^2}\right)^3\Big\}(\lambda_a,0,0)\ne 0$.
Thus, $\frac{\partial^3 g}{\partial \lambda^3}(\lambda_a,0,)\ne
0$. We know that $\frac{\partial g}{
\partial p}(\lambda_a,0,)\ne 0$ since $\frac{\partial^2 F}{\partial t
\partial p}(\lambda_a,a,\dot{q}_a) \ne 0$. It is easy to see that $\frac{\partial^2 F}{\partial t
\partial q}(\lambda_a,a,\dot{q}_a)=0$. Thus $F$ has $A_2^{A_2}$ singularity at $(\lambda_a,a,\dot{q}_a)$ iff
$
\frac{\partial^3 F}{\partial \lambda
\partial q \partial t}(\lambda_a,a,\dot{q}_a)\frac{\partial^3 F}{\partial
t^3}(\lambda_a,a,\dot{q}_a)-\frac{\partial^3 F}{\partial \lambda
\partial t^2}(\lambda_a,a,\dot{q}_a)\frac{\partial^3 F}{\partial q
\partial t^2}(\lambda_a,a,\dot{q}_a)\ne 0
$.
By direct calculation, this is equivalent to
$
\frac{(q_a^+-q_a^-)}{\lambda_a(1-\lambda_a)}\frac{\partial^3
S^+}{\partial (q^+)^3}(q_a^+)\frac{\partial^3 S^-}{\partial
(q^-)^3}(q_a^-)\ne 0
$,
which is not satisfied, since $l(a^+,a^-)$ is $2$-tangent to $L$
at $a^-$.
\end{proof}
\begin{cor}\label{$(1,2m)$-Lagrangian-instability}
Let $L$ be a smooth closed convex curve in
$(\mathbb{R}^{2},\omega)$. The smooth part of $E_{1/2}(L)$ is $(1,2)$-Lagrangian stable, but the cusps of
$E_{1/2}(L)$, seen as part of $GCS(L)$, are not
$(1,2)$-Lagrangian stable; the middle axes and the whole CSS are
not $(1,2)$-Lagrangian stable.
\end{cor}
\begin{rem}
For a convex curve $L\subset\mathbb R^2$,
most singularities which are affine stable are not
affine-Lagrangian stable (compare Theorem \ref{curvestability} and Corollary \ref{$(1,2m)$-Lagrangian-instability}).
Also,
although the cusps  of $E_{1/2}(L)$ are affine-Lagrangian stable
when $E_{1/2}(L)$ is considered by itself,  they are not affine-Lagrangian stable considering  $E_{1/2}(L)\subset GCS(L)$, that is, the meeting of $E_{1/2}(L)$ and  CSS is not affine-Lagrangian
stable.
\end{rem}
\subsection{Discussion}
Because of the large loss of stability for singularities of the GCS,  when going from the affine to the affine-Lagrangian case,  one wonders if it is possible to consider a
coarsen classification of singularities of the GCS of Lagrangian submanifolds, which produces more stable singularities. In fact, the usual Lagrangian equivalence  will do.

As mentioned at the beginning of section \ref{section-sing-GCS}, classification by usual Lagrangian equivalence amounts to considering the unfolding parameters
$y=(\lambda,x)\in\mathbb R\times\mathbb R^{2m}$ on an equal footing.
In this setting, Lagrangian equivalence of $\mathbb E(L)$ and $\mathbb E(\widetilde{L})$ is defined in terms of Lagrangian equivalence of $\mathcal L$ and $\widetilde{\mathcal L}$ in the usual way, which means that their generating families must be
stably $\mathcal R^+$-equivalent (theorem \ref{LagstabR+}), in other words,
there is a symplectomorphism-germ $\Upsilon$ of
$T^*\mathbb R \times T\mathbb R^{2m}$ such that
$\Upsilon(\mathcal L)=\widetilde{\mathcal L}$ and the following
diagram commutes:
$$
\begin{array}{ccccccccc}

& & & Pr &   \\
\mathcal L & \hookrightarrow & T^*\mathbb R \times T\mathbb R^{2m} & \longrightarrow & \mathbb R\times\mathbb R^{2m}
  \\
& & & &  &  \\
& & \downarrow \Upsilon &  & \downarrow     \\
& & & Pr &    \\
\widetilde{\mathcal L} & \hookrightarrow & T^*\mathbb R \times T\mathbb R^{2m} &  \longrightarrow &  \mathbb R\times\mathbb
R^{2m}  \\
\end{array}
$$
where the right-vertical arrow is a diffeomorphism-germ of general form
$$
\mathbb R\times \mathbb R^{2m} \ni
(\lambda,x)\mapsto(\Lambda(\lambda,x),X(\lambda,x)) \in \mathbb R \times
\mathbb R^{2m}.
$$

Comparing with the classifying diagram in Definition \ref{GCS-$(1,2m)$-Lagrangian-eq} for  $(1,2m)$-Lagrangian equivalence, one expects that many singularities of $GCS(L)$
which are Lagrangian stable  are not $(1,2m)$-Lagrangian stable. In fact, for convex Lagrangian curves, it is easy to see that most of the singularities
of Theorem   \ref{curvestability} are Lagrangian stable in the above sense.

However,  the fact that the last projection $\pi: \mathbb R^{1+2m} \to \mathbb R^{2m}$ is not taken into account is an obvious  indication that usual Lagrangian equivalence is not the correct equivalence relation  for classification of the singularities of $GCS(L)$,  because this latter is the image under $\pi$ of the locus of critical points of $\pi$ restricted to $\mathbb E(L)$.

This becomes even clearer when we also analyze the non-symplectic case. In this case, consider the following {\it extended chord transformation}
$$\Gamma: \mathbb R \times\mathbb R^n\times \mathbb R^n\to\mathbb R \times T\mathbb R^n \ , \ (\lambda,x^+,x^-)\mapsto (\lambda,\Gamma_{\lambda}(x^+,x^-)) \ ,$$
where $\Gamma_{\lambda}:\mathbb R^n\times \mathbb R^n\to T\mathbb R^n $ is a simpler $\lambda$-chord transformation,
\begin{equation}\label{standardchordtransf}
\Gamma_{\lambda}(x^+,x^-)=(x,\dot x)=\left(\lambda x^+ + (1-\lambda)x^-,\frac{x^+-x^-}{2}\right) \ ,
\end{equation}
which  differs from $\Phi_{\lambda}$ only in the kind of linear equation for $\dot x$, compare (\ref{standardchordtransf}) to  (\ref{x}) and (\ref{genx.}), this latter chosen in the symplectic case so that $(\Phi_{\lambda}^{-1})^*(\delta_{\lambda}\omega)=\dot\omega$ (no extra semi-basic form in the r.h.s.).

Now, let $M$ and $\widetilde{M}$ be germs of $m$-dimensional smooth submanifolds of $\mathbb R^n$, $n\leq 2m$, and let $\mathbb M$ and $\widetilde{\mathbb M}$ be the chord transformed cylinders
$$\mathbb M=\Gamma(\mathbb R\times M\times M) \  , \  \widetilde{\mathbb M}=\Gamma(\mathbb R\times\widetilde{M}\times\widetilde{M}) \ .$$

\begin{defn}\label{eq-gcs} Germs of $GCS(M)$ and
$GCS(\widetilde{M})$ are {\bf chord equivalent} if there is a diffeomorphism-germ $\Theta$ of $\mathbb R\times T\mathbb R^n$ s.t. $\widetilde{\mathbb M}=\Theta(\mathbb M)$ and the
following commutes:
$$
\begin{array}{ccccccc}
 & id_{\mathbb R}\times pr & & \pi & \\
 \mathbb R\times
T\mathbb R^n & \longrightarrow & \mathbb R \times \mathbb R^n &
\longrightarrow & \mathbb
 R^n \\
 & & & & & & \\
 \downarrow \Theta &  & \downarrow &  & \downarrow  \\
 & id_{\mathbb R}\times pr & & \pi & \\
 \mathbb R\times T\mathbb R^n &
\longrightarrow & \mathbb R \times \mathbb R^n & \longrightarrow &
\mathbb
 R^n \\
\end{array}
$$
where \emph{vertical} arrows indicate diffeomorphism-germs, as follows:
$$
\Theta: \mathbb R\times T\mathbb R^n \ni
(\lambda,x,\dot{x})\mapsto(\Lambda(\lambda,x),X(x),\dot{X}(\lambda,x,\dot{x}))
\in \mathbb R \times T\mathbb R^n,
$$
$$
\mathbb R\times \mathbb R^n \ni
(\lambda,x)\mapsto(\Lambda(\lambda,x),X(x)) \in \mathbb R \times
\mathbb R^n,
$$
$$
\mathbb R^n \ni x\mapsto X(x) \in \mathbb R^n.
$$
\end{defn}
\begin{defn}
A singularity of $GCS(M)$ is {\bf affine stable} if it is a stable singularity  under its classification by chord equivalence.
\end{defn}

Given the above definition, by somewhat lengthy but straightforward computations, one proves Theorem  \ref{curvestability} for the $GCS$ of convex curves using classification by chord equivalence. The classification of the singularities of $GCS(M)$ in the other known cases, for instance hyperplanes, can be similarly accomplished by chord equivalence, which gives the correct affine-invariant classification of the singularities of $GCS(M)$ for general $m$-dimensional submanifolds $M\subset\mathbb R^n$, $n\leq 2m$.

Comparison of the classifying diagram in Definition \ref{eq-gcs} for chord equivalence with the classifying diagram in Definition \ref{GCS-$(1,2m)$-Lagrangian-eq} for  $(1,2m)$-Lagrangian equivalence shows their obvious analogy.
On the other hand, the obvious analog of the classifying diagram for usual Lagrangian equivalence, when no symplectic form has to be accounted for, is
$$
\begin{array}{ccccccc}
 & id_{\mathbb R}\times pr &  \\
 \mathbb R\times
T\mathbb R^n & \longrightarrow & \mathbb R \times \mathbb R^n  \\
 & & & & & & \\
 \downarrow \Theta &  & \downarrow   \\
 & id_{\mathbb R}\times pr &  \\
 \mathbb R\times T\mathbb R^n &
\longrightarrow & \mathbb R \times \mathbb R^n  \\
\end{array}
$$
where vertical arrows indicate diffeomorphism-germs of the form:
$$
\Theta: \mathbb R\times T\mathbb R^n \ni
(\lambda,x,\dot{x})\mapsto(\Lambda(\lambda,x),X(\lambda,x),\dot{X}(\lambda,x,\dot{x}))
\in \mathbb R \times T\mathbb R^n,
$$
$$
\mathbb R\times \mathbb R^n \ni
(\lambda,x)\mapsto(\Lambda(\lambda,x),X(\lambda,x)) \in \mathbb R \times
\mathbb R^n.
$$

Of course, applying the above {\it wrong} equivalence relation to classify singularities of $GCS(M)$ for general $m$-dimensional submanifolds $M\subset\mathbb R^n$, $n\leq 2m$, produces many more stable singularities than when applying the {\it correct} classifying diagram in Definition  \ref{eq-gcs}.

Thus, choosing the correct classifying diagram in both the non-symplectic and the symplectic cases shows that most singularities of the $GCS$ which are stable when no symplectic form has to be accounted for, cease to be stable when there is a symplectic form to be accounted for. In other words, there is breakdown of stability due to a symplectic form. Other similar cases, of breakdown of simplicity due to a symplectic form, can
be found in \cite{D,DJZ} and specially in \cite{DR}.

\end{document}